\documentclass[a4paper,12pt,twoside]{amsart}

\usepackage[cm]{fullpage}
\usepackage{mymacros,oitp}
\usepackage[%dvipdfm,%
pagebackref]{hyperref}

\includecomment{comment}

\usepackage{comment} 
\title[Simple singularities]{Homological mirror symmetry
for Milnor fibers of simple singularities}
\author[Y.~Lekili]{Yank\i\ Lekili}
\address{
Department of Mathematics
Imperial College London
South Kensington 
London
SW7 2AZ}
\email{y.lekili@imperial.ac.uk}
\author[K.~Ueda]{Kazushi Ueda}
\address{
Graduate School of Mathematical Sciences,
The University of Tokyo,
3-8-1 Komaba,
Meguro-ku,
Tokyo,
153-8914,
Japan.}
\email{kazushi@ms.u-tokyo.ac.jp}
\date{}

\begin{document}

\begin{abstract}
We prove homological mirror symmetry
for Milnor fibers of simple singularities
in dimensions greater than one,
which are among the log Fano cases of
\cite[Conjecture 1.5]{1806.04345}.
The proof is based on a relation
between matrix factorizations and Calabi--Yau completions.
As an application,
we give an explicit computation of the Hochschild cohomology group
of the derived $n$-preprojective algebra
of a Dynkin quiver for any $n \ge 1$,
and the symplectic cohomology group
of the Milnor fiber of any simple singularity in any dimension
greater than one.
\end{abstract}

\maketitle

\section{Introduction}

%A Kleinian singularity
%is a surface singularity
%appearing in the quotient of $\bC^2$
%by a finite subgroup of $\SL_2(\bC)$,
%and known to have many characterizations
%such as a rational double point,
%a rational Gorenstein singularity,
%a canonical singularity,
%or a simple singularity.
%Finite subgroups of $\SL_2(\bC)$ are also known as
%binary polyhedral groups,
%and have the well-known ADE classification.
%The coordinate ring of the quotient
%%(i.e., the invariant ring)
%is generated by three elements,
%which can be chosen in such a way that the relation is given by

A \emph{simple singularity} is an isolated hypersurface singularity of modality zero.
Arnold classified such singularities;
up to right equivalence,
they are given by one of the following:
\begin{align} \label{eq:simple}
\begin{split}
 A_\ell &\colon x_1^{\ell+1} + x_2^2 + \cdots + x_{n+1}^2 = 0, \quad \ell = 1, 2, \ldots \\
%  \label{eq:Aell} \\
 D_\ell &\colon x_1^{\ell-1} + x_1 x_2^2 + x_3^2+ \cdots + x_{n+1}^2 = 0, \quad \ell = 4, 5, \ldots \\
%  \label{eq:Dell} \\
 E_6 &\colon x_1^4 + x_2^3 + x_3^2 + \cdots + x_{n+1}^2 = 0, \\
%  \label{eq:E6} \\
 E_7 &\colon x_1^3 + x_1 x_2^3 + x_3^2 + \cdots + x_{n+1}^2 = 0, \\
%  \label{eq:E7} \\
 E_8 &\colon x_1^5 + x_2^3 + x_3^2 + \cdots + x_{n+1}^2 = 0.
%  \label{eq:E8}
\end{split}
\end{align}
In the case $n=2$,
simple surface singularities have many other characterizations,
such as Kleinian singularities,
rational double points,
or canonical singularities,
to name a few.

Let $\wv$ be one of these defining polynomials,
which we think of as a holomorphic function on $\bC^{n+1}$,
and equip $\wv^{-1}(1)$ with the Liouville structure
induced from the standard one on $\bC^{n+1}$.
This is the Liouville completion of the \emph{Milnor fiber},
which is the Liouville domain obtained by intersecting $\wv^{-1}(1)$ with a ball.
Let $\cW \lb \wv^{-1}(1) \rb$ denote
the idempotent-complete derived wrapped Fukaya category
of $\wv^{-1}(1)$.

For $n \geq 2$, since $\wv^{-1}(1)$ is not a log Calabi--Yau manifold
but a log Fano manifold,
its mirror is not a manifold
but a \emph{Landau--Ginzburg model},
by which we mean a pair of a stack and a section of a line bundle on it.
One way to obtain a Landau--Ginzburg mirror
of a log Fano manifold
is to first remove a divisor to make it log Calabi--Yau,
then find its mirror, which is another log Calabi--Yau manifold,
and finally add a potential to this mirror
\cite{MR2386535,MR2537081}.
This produces a Landau--Ginzburg mirror
whose underlying manifold is of the same dimension
as the original manifold.
%In the case of surfaces,
%the class of log Calabi--Yau varieties
%is the same as that of cluster varieties,
%and mirror symmetry for them is studied
%in depth (see \cite{MR3415066}).
When the singularity is toric
(i.e., a simple surface singularity of type A),
there is a standard choice
for the divisor to remove,
and the resulting mirror is the Landau--Ginzburg model
consisting of a complement of a toric divisor
in the minimal resolution of the singularity of the same type
and a monomial function on it
(see e.g.~\cite[Section 9.2]{MR3502098}).
%For a general Kleinian singularity
%(including the type A case),
The choice of the divisor is not unique in general,
and there are multiple mirrors for a given Milnor fiber.

In this paper,
we consider an alternative mirror
of the Milnor fiber of a simple singularity
based on transposition of invertible polynomials
introduced in \cite{MR1214325,MR1310310}.
A weighted homogeneous polynomial
$
 \w \in \bC[x_1,\ldots,x_{n+1}]
$
with an isolated critical point at the origin
is \emph{invertible}
if there is an integer matrix
$
 A = (a_{ij})_{i, j=1}^{n+1}
$
with non-zero determinant
such that
\begin{align}
 \w = \sum_{i=1}^{n+1} \prod_{j=1}^{n+1} x_j^{a_{ij}}.
%  \in \bC[x_1,\ldots,x_n]
\end{align}
The \emph{transpose} of $\w$ is defined as
\begin{align}
 \wv =  \sum_{i=1}^{n+1} \prod_{j=1}^{n+1} x_j^{a_{ji}},
\end{align}
whose exponent matrix $\Av$
is the transpose matrix of $A$.
The group
\begin{align}
 \Gamma_{\w} \coloneqq
 \lc (t_0, t_1, \ldots, t_{n+1}) \in (\Gm)^{n+2} \relmid
  t_1^{a_{1,1}} \cdots t_{n+1}^{a_{1,{n+1}}}
   = \cdots
   = t_1^{a_{{n+1},1}} \cdots t_{n+1}^{a_{{n+1},{n+1}}}
   = t_0 t_1 \cdots t_{n+1}
   \rc \end{align}
acts naturally on
$
 \bA^{n+2} \coloneqq \Spec \bC[x_0,\ldots,x_{n+1}].
$
Let
$
 \mf \lb \bA^{n+2}, \Gamma_{\w}, \w + x_0 \cdots x_{n+1} \rb
$
denote the idempotent completion
of the dg category of $\Gamma_{\w}$-equivariant coherent matrix factorizations of
$\w + x_0 \cdots x_{n+1}$
on $\bA^{n+2}$
in the sense of \cite{MR3366002}.
\pref{cj:hms} below is given
in \cite[Conjecture 1.5]{1806.04345}:

\begin{conjecture} \label{cj:hms}
For any invertible polynomial $\w$,
one has a quasi-equivalence
\begin{align} \label{eq:hms}
 \mf \lb \bA^{n+2}, \Gamma_\w, \w+x_0 \cdots x_{n+1} \rb
  \simeq \cW \lb \wv^{-1}(1) \rb.
\end{align}
\end{conjecture}

In other words,
the Landau--Ginzburg model
$
 \lb \ld \bA^{n+2}/\Gamma_{\w} \rd, \w + x_0 \cdots x_{n+1} \rb
$
is mirror to the Liouville manifold $\wv^{-1}(1)$.
The main result of this paper is the following:

\begin{theorem} \label{th:main}
\pref{cj:hms} holds for $n \ge 2$ and $\wv$ one of the defining polynomials
of simple singularities
appearing in \pref{eq:simple}.
\end{theorem}

The proof of \pref{th:main}
consists of four steps.
The first step is the quasi-equivalence
\begin{align} \label{eq:removing xyz}
  \mf(\bA^{n+2}, \Gamma_{\w}, \w+x_0 \cdots x_{n+1})
  \simeq
  \mf(\bA^{n+2}, \Gamma_{\w}, \w),
\end{align}
which comes from the fact that $\w+x_0 \cdots x_{n+1}$
is right equivalent to $\w$
by a formal change of variables,
which holds if $n \ge 2$
and $\w$ defines a simple singularity.

The second step is the quasi-equivalence
\begin{align} \label{eq:Pi mf}
  \mf(\bA^{n+2}, \Gamma_{\w}, \w)
  \simeq
  \Pi_n \lb \mf(\bA^{n+1}, \Gamma_{\w}, \w) \rb,
\end{align}
where $\Pi_n$ denotes the $n$-Calabi--Yau completion
in the sense of \cite{MR2795754}.
This holds for any invertible polynomial $\w$
and any $n \ge 0$.

The third step is the quasi-equivalence
\begin{align} \label{eq:mf=perf AQ}
  \mf(\bA^{n+1}, \Gamma_{\w}, \w)
  \simeq
  \operatorname{perf} A_Q
\end{align}
with
the dg category
$\perf A_Q$
of perfect dg modules
over the path algebra $A_Q$
of a Dynkin quiver $Q$
(with any orientation)
of the corresponding type.
For type A,
this is proved in
\cite[Theorem 3.1]{0506347}
for $n=0$,
and the $n \ge 1$ case
follows
either from the $n=0$ case
and the Kn\"orrer periodicity
\cite{MR877010},
or as a special case of
\cite[Theorem 1.2]{MR2803848}.
For type D,
this follows from \cite[Section 4]{MR3030671}.
For type E,
this follows either from the combination of
\cite[Theorem 1]{1903.01351}
and
\cite[Proposition 3.4]{MR1882336}
or by finding a suitable mutation of a generator
appearing in \cite[Theorem 2]{1903.01351}.
Note that
\cite[Theorem 3.1]{MR2313537}
gives a result
close to \pref{eq:mf=perf AQ},
which is not exactly the same
since the grading group is different.

The last step is
\begin{align} \label{eq:W=AQ}
  \cW \lb \wv^{-1}(1) \rb
  \simeq
  \Pi_n \lb \perf A_Q \rb,
\end{align}
which holds if $n \ge 2$ and $\wv$ defines a simple singularity.
As discussed in \pref{sc:wrap},
the proof of \pref{eq:W=AQ} in \cite{MR3692968}
depends on the computation
of the Hochschild cohomology of $\Pi_n \lb \perf A_Q \rb$,
which was missing for type E cases in \cite{MR3692968}
and is done in \pref{sc:HH}.

For $n=1$
not covered by \pref{th:main},
a quasi-equivalence of
the full subcategory $\cF(\wv^{-1}(0))$
of $\cW(\wv^{-1}(0))$
consisting of
(direct summands of bounded complexes of)
compact Lagrangians
and a category $\perf Z_{\w}$
equivalent to the full subcategory of 
$
\mf \lb \bA^{3}, \Gamma_\w, \w+x_0 x_1 x_2 \rb
$
consisting of homologically finite objects
(i.e., those $X$ satisfying $\dim \bigoplus_{i \in \bZ} \Ext^i(X, Y) < \infty$
for any object $Y$)
is given in \cite[Theorem 1.1]{2003.01106}.

As an application of \pref{eq:Pi mf},
we compute the Hochschild cohomology group
of the $n$-Calabi--Yau completion
$\Pi_n(A_Q)$,
also known as the \emph{derived $n$-preprojective algebra},
of the path algebra $A_Q$
of any Dynkin quiver $Q$
for any $n \ge 1$.
It is possible to compute the Hochschild homology
along the same line.

The zero-th cohomology of the derived $2$-preprojective algebra
is the \emph{preprojective algebra}.
The Hochschild homology and cohomology
of the preprojective algebra
of the path algebra of a Dynkin quiver
is calculated in
\cite{MR1648626,MR1632808,MR2372205}.
Even the calculus structure
in the sense of \cite{MR1039918,MR1783778}
(which includes the Batalin--Vilkovisky structure
and is known to be derived invariant
\cite{MR3945161})
is calculated in \cite{MR2561764},
and it is an interesting problem to do the same
for the derived $n$-preprojective algebra.
Note that
the preprojective algebra and the derived $2$-preprojective algebra
of a Dynkin quiver are very different.
The derived $2$-preprojective algebra
of a Dynkin quiver is a smooth dg algebra,
which has cohomology in every negative cohomological degree,
and moreover is not formal.
In contrast,
the preprojective algebra is always concentrated in cohomological degree $0$ by definition,
and the global dimension
is infinite for a Dynkin quiver.

It follows from
\cite[Theorem 1.1]{MR3121862},
combined with
\cite[Theorem 1.4]{1712.09126}
which builds on \cite{MR3121862, 1712.00225},
that
the closed-open map
of any Weinstein manifold
from the symplectic cohomology
to the Hochschild cohomology of the wrapped Fukaya category
is an isomorphism:
\begin{align} \label{eq:SH_HH}
 \operatorname{SH}^*(M) \simto \HH^*(\cW(M)).
\end{align}
Hence, by Theorem \ref{th:main}, we see that the symplectic cohomology
of the Milnor fiber
$\wv^{-1}(1)$
of a simple singularity
for $n \ge 2$
is isomorphic to
$
 \HH^* \lb \Pi_n \lb A_Q \rb \rb.
$
This enables us to give an explicit computation of the symplectic cohomology of Milnor fibers
of all simple singularities in a uniform way.
Previous partial results computing symplectic cohomology
for Milnor fibers of simple simple singularities appeared
in \cite{MR3692968} for $A_\ell$ and $D_\ell$ in complex dimension 2,
and in \cite{MR3483060, MR3489066},
for various versions of symplectic cohomology
for certain higher dimensional $A_\ell$-Milnor fibers
for which an associated Morse--Bott spectral sequence yields computations.
Our computation also shows that
$\HH^* \lb \Pi_1 \lb A_Q \rb \rb$
is not isomorphic to $\operatorname{SH}^*(\wv^{-1}(0))$
given in \cite[Section 3.3]{2003.01106},
which is consistent with 
the failure of \pref{eq:W=AQ} for $n=1$.

This paper is organized as follows:
In \pref{sc:quiver},
we collect basic definitions and results on Calabi--Yau completions
and trivial extension algebras.
In \pref{sc:wrap},
we recall the description
of the wrapped Fukaya category
of the Milnor fiber of a simple singularity
for $n \ge 2$
in terms of the $n$-Calabi--Yau completion
of a Dynkin quiver
of the corresponding type.
In \pref{sc:CY-completion},
we prove \pref{eq:removing xyz} and \pref{eq:Pi mf}.
The computation of Hochschild cohomologies
of the derived preprojective algebras
of Dynkin quivers
are given
in \pref{sc:HH}.

\textit{Acknowledgment}:
We thank the anonymous referees for reading the manuscript carefully
and suggesting many improvements and corrections.
Y.~L.~is partially supported by the Royal Society URF\textbackslash R\textbackslash180024. K.~U.~is partially supported by Grant-in-Aid for Scientific Research
(15KT0105, 16K13743, 16H03930).

\section{Calabi--Yau completions and trivial extension algebras}
 \label{sc:quiver}

%A \emph{quiver} $Q = (Q_0, Q_1, s, t)$
%consists of
%a set $Q_0$ of \emph{vertices},
%a set $Q_1$ of \emph{arrows},
%and two maps $s, t \colon Q_1 \to Q_0$
%sending an arrow to its \emph{source} and \emph{target} respectively.
%The \emph{path algebra} of $Q$ will be denoted by $A_Q$.

The \emph{$n$-Calabi--Yau completion}
(or the \emph{derived $n$-preprojective algebra})
of a dg category
$
 \scrA
$
is defined in \cite[Section 4.1]{MR2795754}
as the tensor algebra
\begin{align}
 \Pi_{n} (\scrA)
  &\coloneqq T_{\scrA}(\theta)
  \coloneqq \scrA \oplus \theta \oplus \theta \otimes_\scrA \theta \oplus \cdots,
\end{align}
where
the $\scrA$-bimodule
$
 \theta
  \coloneqq \Theta[n-1]
$
is a shift of the inverse dualizing complex
$
 \Theta
  \coloneqq \hom_{\scrA^{\mathrm{e}}} (\scrA, \scrA^{\mathrm{e}}).
%  \simeq \bS^{-1}
$

A dg algebra is regarded as a dg category with one object.
The Morita invariance of the Calabi--Yau completion
shown in \cite[Proposition 4.2]{MR2795754}
implies that Calabi--Yau completion commutes
with the operation of
taking the dg category of perfect dg modules:
\begin{align}
  \Pi_n \, (\perf \scrA) \simeq \perf (\Pi_n\, \scrA).
\end{align}

The \emph{Ginzburg dg algebra}
$
 \scrG_Q^{n}
$
of a quiver $Q$
(without potential)
is a model of the $n$-Calabi--Yau completion
$\Pi_{n} \lb A_Q \rb$
of the path algebra $A_Q$,
defined in
\cite[Section 6.2]{MR2795754}
after \cite{0612139}
as the path algebra of the graded quiver
$\overline{Q}$
with same vertices as $Q$ and arrows consisting of
\begin{itemize}
 \item
the original arrows $g \in Q_1$ in degree $1$,
 \item
the opposite arrows $g^*$ for each arrow $g \in Q_1$
in degree $1-n$, and
 \item
loops $h_v$ at each vertex $v \in Q_0$
in degree $1-n$,
\end{itemize}
equipped with the differential $d$
given by
\begin{align}
 dg = dg^* = 0  \ \mbox{ and } \  dh = \sum_{g\in Q_1} g^* g - g g^*
\end{align}
where $h = \sum_{v\in Q_0} h_v$.

The \emph{degree $n$ trivial extension algebra}
of a finite-dimensional algebra $A$
is defined as
$
% A^{n-1} \coloneqq 
 A \oplus A^\dual[-n]
$
equipped with the multiplication
$
 (a, f) \cdot (b, g) = (ab, ag+fb),
$
where
$
 A^\dual
$
is the dual of $A$ as a vector space.

The degree $n$ trivial extension algebra
$
 B_Q^{n}
$
of the path algebra $A_Q$ of a Dynkin quiver $Q$
is the (derived) Koszul dual of $\scrG^{n}_Q$
in the sense that
\begin{align} \label{eq:Koszul_duality}
 \hom_{\scrG_Q^{n}} \lb \bfk_{\scrG}, \bfk_{\scrG} \rb
  \simeq B_Q^{n},
 \qquad
        \hom_{\left(B_Q^{n}\right)^\mathrm{op}} \lb \bfk_B, \bfk_B \rb
  \simeq \left( \scrG_Q^{n} \right)^\mathrm{op},
\end{align}
where
$
 \bfk_{\scrG} \coloneqq \bigoplus_{v \in Q_0} S_v
$
is the direct sum of simple left $\scrG_Q^{n}$-modules $S_v$
associated with vertices $v \in Q_0$,
and similarly for $\bfk_B$
(see e.g.~\cite[Theorem 23, Corollary 25]{MR3692968}).

This Koszul duality implies an isomorphism
\begin{align} \label{eq:hh}
 \HH^* \lb \scrG_Q^{n} \rb
  \cong \HH^* \lb B_{Q}^{n} \rb
\end{align}
of Hochschild cohomologies
(see e.g.~\cite[Theorem 1]{MR2134291} and \cite[Theorem 3.4]{MR3941473}).

\section{Wrapped Fukaya category of the Milnor fiber of simple singularity}
\label{sc:wrap}

Let $\wv$ be one of the defining polynomials of a simple singularity and $M^{n} = \wv^{-1}(1)$ be the Milnor fiber, which we view as a Weinstein manifold where the Weinstein structure is induced by restriction from the ambient $\bC^{n+1}$. It is well known that this Weinstein manifold is symplectomorphic (in fact, Weinstein homotopic) to the plumbing $X_Q$ of cotangent bundles of spheres $T^*S^{n}$ according to the Dynkin diagram $Q$ corresponding to the simple singularity. One way to see this is to verify it directly for $n=1$, and then use the fact that in higher dimensions the Milnor fiber is obtained by stabilization --- increasing the dimension corresponds to suspension of the Lefschetz fibration \cite{MR2651908}. See also \cite{MR2786590} for an explicit construction of a symplectic structure on plumbings. This stabilization point of view also enables one to describe $M$ via Legendrian surgery. Namely $M$ is obtained by attaching critical handles to a Legendrian link $\Lambda_Q^{n-1}$ on $\partial \mathbb{D}^{n}$ whose components are unknotted Legendrian spheres $S^{n-1}$ which are clasped together (as in Hopf link) according to the Dynkin diagram $Q$. The direct sum of co-cores to the critical handles (i.e., cotangent fibers away from the plumbing region) form a generating object of the wrapped Fukaya category by the main theorem in \cite{1712.09126}, and the surgery formula of \cite{MR2916289, 1906.07228} allows one to explicitly compute the endomorphism algebra of this generator as the Chekanov--Eliashberg algebra $\CE^*(\Lambda_Q^{n-1})$. 

This Chekanov--Eliashberg algebra was computed directly in the case $n=2$ in the paper \cite{MR3692968} and the resulting dg algebra was shown to be quasi-isomorphic to the derived multiplicative preprojective algebra of the corresponding Dynkin type.
Moreover, working over $\bC$, for $Q= A_\ell$ or $D_\ell$,
it was shown in \cite[Theorem 13]{MR3692968} that the derived multiplicative preprojective algebra of Dynkin type $Q$ is quasi-isomorphic to the Ginzburg algebra $\mathscr{G}^2_Q$,
also known as the derived (additive) preprojective algebra of Dynkin type $Q$.
It was conjectured in op.~cit.~that the same result holds for $Q=E_6, E_7, E_8$ and this is indeed so. The key ingredient for the proof of \cite[Theorem 13]{MR3692968} to go through that was missing in the case $Q=E_6,E_7,E_8$ was the computation that
\begin{align}
  \HH^2 \lb \mathscr{G}^2_Q\rb^s =0 \text{ for } s <0,
\end{align}
but this follows from computations given in Section \ref{sc:HH} below. 

For $n \ge 3$, one can do a direct computation in an analogous way, but we can also deduce this by the Koszul duality result given in \cite[Theorem 58]{1701.01284} which shows that $\CE^*(\Lambda_Q^{n-1})$ is the (derived) Koszul dual of the endomorphism algebra of the union of the core spheres of the plumbing. Notice that for $n \ge 3$, $\wv$ is suspended at least twice, thus the formality of the endomorphism algebra of vanishing cycles in the compact Fukaya category of $\wv^{-1}(1)$ follows automatically by \cite[Proposition 4.4]{MR2651908}
(the formality of the $A_\infty$-algebra $\cA$ and a $\cA$-bimodule $\cB/\cA$ in Seidel's notation
is obvious in the case at hand,
since $\Gamma$ is a tree and one can shift the objects
to put all morphisms in degree $0$).
Putting it all together, we conclude that
$\CE^*(\Lambda_Q^{n-1})$ is Koszul dual
to the degree $n$ trivial extension algebra $B_Q^{n}$
of the path algebra $A_Q$
of a Dynkin quiver of the corresponding type
(see also \cite{MR3977875} for another example).

As a result of these computations,
for $n \ge 2$
we have a quasi-isomorphism
\begin{align} \label{eq:iso}
 \CE^* \lb \Lambda_Q^{n-1} \rb \simeq \scrG^{n}_Q
\end{align}
over $\bC$,
which implies a quasi-equivalence
\begin{align} \label{eq:W_Pi}
 \cW(\wv^{-1}(1)) \simeq \perf \Pi_{n}(A_Q)
\end{align}
between the wrapped Fukaya category of $\wv^{-1}(1)$
and the dg category of perfect modules over $\Pi_{n}(A_Q)$.

\begin{remark}
Note
from \cite[Proposition 3.4]{MR1882336}
that $A_Q$ is derived equivalent to the Fukaya--Seidel category
$\cF(\wv)$
of the LG-model
$
 \wv \colon \bC^{n+1} \to \bC.
$
Thus \pref{eq:W_Pi} shows that
$\cW \lb \wv^{-1}(1) \rb$
is the Calabi--Yau completion of $\cF(\wv)$
for $n \ge 2$.
Although this relationship between $\cF(\wv)$ and $\cW \lb \wv^{-1}(1) \rb$
is not true in general,
we expect it to hold
when $\wv$ is a double suspension of an invertible polynomial
whose Milnor fiber is a log Fano manifold,
since one has
\begin{multline}
  \w(x_1,\ldots,x_{n-1}) + x_n^2 + x_{n+1}^2 + x_0 \cdots x_{n+1} \\
  =\w(x_1,\ldots,x_{n-1})
  +\left( \sqrt{1-\frac{1}{4} (x_0 \cdots x_{n-1})^2} x_n \right)^2
  +\left( x_{n+1}+\frac{1}{2} x_0 \cdots x_n \right)^2
\end{multline}
in
$
\bfk \db[ x_0,\ldots, x_{n+1} \db].
$
\end{remark}

\begin{remark}
The isomorphism \pref{eq:iso} remains true for $n \ge 3$ over an arbitrary commutative ring,
but for $n=2$ we have to require that $2$ is invertible for type $D_\ell, E_6,E_7,E_8$,
$3$ is invertible for type $E_6,E_7, E_8$, and
5 is invertible for type $E_8$.
Otherwise, $\CE^*(\Lambda_Q)$ is quasi-isomorphic to the derived multiplicative preprojective algebra (see \cite{MR4033516}) which is not quasi-isomorphic to the derived (additive) preprojective algebra $\Pi_{n}(A_Q)$. 
\end{remark}

\section{Matrix factorizations and Calabi--Yau completions}
 \label{sc:CY-completion}

Let $\Gamma$ be a subgroup
of $(\Gm)^{n+1}$
acting diagonally on
$
 \bA^{n+1} \coloneqq \Spec \bC[x_1,\ldots,x_{n+1}].
$
Assume that $\Gamma$ is a finite extension
of the multiplicative group $\Gm$,
so that the group 
$
 \Char(\Gamma)
  \coloneqq \Hom(\Gamma, \Gm)
$
of characters of $\Gamma$
is an extension of a finite group by $\bZ$.
The coordinate ring $\bC[x_1,\ldots,x_{n+1}]$
has a $\Char(\Gamma)$-grading
coming from the $\Gamma$-action on $\bA^{n+1}$,
and we set $\chi_i \coloneqq \deg x_i$
for $i \in \{ 1, \ldots, n+1 \}$.
Let
$
 \w \in \bC[x_1,\ldots,x_{n+1}]_\chi
  % \coloneqq \lb \bC[x_1,\ldots,x_{n+1}] \otimes \chi \rb^\Gamma
$
be a homogeneous element
of degree
$
 \chi
  \in
   \Char(\Gamma).
$
%Let $\kappa$ be the character
%obtained as the determinant of the $\Gamma$-action
%on $\bC x_1 \oplus \cdots \oplus \bC x_{n+1}$.
Assume that $\w$ has an isolated critical point at the origin,
so that
the structure sheaf $\cO_0$ of the origin
split-generates
$
 \mf \lb \bA^{n+1}, \w \rb
$
%the residue field $\bC[x_1,\ldots,x_{n+1}]/(x_1,\ldots,x_{n+1})$
%by shifts, cones, and direct summands
by \cite[Proposition A.2]{MR2776613}
(see also \cite{MR2735755,MR2824483}).
Let
$
 R \subset \Char(\Gamma)
$
be a set of representatives of the group
$
 \Char(\Gamma)/(\chi),
$
which we assume to be finite.
Then
%the direct sum
$
 \cE \coloneqq \bigoplus_{\rho \in R} \cO_0(\rho)
$
split-generates
$
 \mf \lb \bA^{n+1}, \Gamma, \w \rb,
$
since the autoequivalence
$
 M \mapsto M(\chi)
$
of
$
 \mf \lb \bA^{n+1}, \Gamma, \w \rb
$
shifting the $\Gamma$-weight by $\chi$
is isomorphic to the functor $M \mapsto M[2]$
shifting the cohomological grading by 2.

The $n$-Calabi--Yau completion
of the dg Yoneda algebra
$
 \scrA \coloneqq \hom(\cE, \cE)
$
is given by
\begin{align}
 \Pi_{n} (\scrA)
%  &\coloneqq T_{\scrA}(\theta)
  \coloneqq \scrA \oplus \theta \oplus \theta \otimes_\scrA \theta \oplus \cdots
  \simeq \bigoplus_{i=0}^\infty \hom(\cE, \theta^i(\cE))
\end{align}
where
$\theta = \Theta[n-1]$ as in \pref{sc:quiver},
and we abuse notation
and use the same symbol for an autoequivalence and its graph bimodule.
Since $\Theta$ is the graph
of the inverse Serre functor $\bS^{-1}$,
we have
\begin{align} \label{eq:theta_S}
 \theta = \bS^{-1}[n-1].
\end{align}

Now, as in \cite[Section 2]{1806.04345},
we introduce another variable $x_0$
of degree
$
 \chi_0
  \coloneqq \chi - (\chi_1 + \cdots + \chi_{n+1}),
%  \in \Char(\Gamma)
$
and
consider the polynomial ring
$
 \bC[x_0,x_1,\ldots,x_{n+1}]
$
in $n+2$ variables,
which naturally contains
$
 \bC[x_1,\ldots,x_{n+1}]
$
as a subring.
One has
\begin{align} \label{eq:kunneth}
  \mf(\bA^{n+2},\w) \simeq \mf(\bA^1,0) \otimes \mf(\bA^{n+1},\w)
\end{align}
e.g., by the ungraded ($G=H=1$) version of
\cite[Lemma 3.52]{MR3270588}
with $v=0$;
note that
$\mf(\bA^1,0)$ is obtained from $\coh \bA^1$
by collapsing the cohomological grading to $\bZ/2 \bZ$,
and the tensor product
of split-generators of
$\mf(\bA^1,0)$ and $\mf(\bA^{n+1},\w)$
gives a split-generator of $\mf(\bA^{n+2},\w)$
since the critical locus of $\w$
as a function on $\bA^{n+2}$
is the product
of $\bA^1$ times that as a function on $\bA^{n+1}$.

As shown in \cite[Theorem 2.5]{MR3063907}
whose proof carries over directly to $\Gamma$-graded cases,
graded Auslander--Reiten duality
\cite{MR915178} implies that
\begin{align} \label{eq:mf_Serre}
 \bS \coloneqq (\chi_0)[n-1]
\end{align}
is a Serre functor on $\mf \lb \bA^{n+1}, \Gamma, \w\rb$.
It follows from \pref{eq:theta_S} and \pref{eq:mf_Serre}
that
\begin{align} \label{eq:theta_chi}
 \theta \simeq (-\chi_0).
\end{align}

Let $\cF$ be the generator of $\mf \lb \bA^{n+2}, \Gamma, \w \rb$
obtained from the tensor product of the generator $\cE$ of
$
 \mf \lb \bA^{n+1}, \Gamma, \w \rb
$
and the generator $\bC[x_0]$ of $\coh \bA^1$.
If we write both of the forgetful functors
$
 \mf \lb \bA^{n+1}, \Gamma, \w \rb \to \mf \lb \bA^{n+1}, \w \rb
$
and
$
 \mf \lb \bA^{n+2}, \Gamma, \w \rb \to \mf \lb \bA^{n+2}, \w \rb
$
as $\overline{(\bullet)}$,
%by a slight abuse of notation,
then one has
\begin{align} \label{eq:non-equiv}
 \hom \lb \overline{\cF},\overline{\cF} \rb
  \simeq \hom \lb \overline{\cE}, \overline{\cE} \rb \otimes \bC[x_0]
  \simeq \bigoplus_{\rho \in \Char(\Gamma)} \hom(\cE, \cE(\rho)) \otimes \bC[x_0].
\end{align}
%through \pref{eq:kunneth}.
Since $\deg (x_0) = \chi_0$,
by taking the $\Gamma$-invariant part of \pref{eq:non-equiv}
and using \pref{eq:theta_chi},
one obtains
\begin{align}
 \hom(\cF,\cF)
  \simeq \bigoplus_{i=0}^\infty \hom(\cE, \cE(-i \chi_0))
  \simeq \bigoplus_{i=0}^\infty \hom(\cE, \theta^i(\cE))
  \simeq \Pi_{n}(\scrA),
\end{align}
which shows
the quasi-equivalence \pref{eq:Pi mf}.

If $n$ is greater than one,
then the degree of $x_1 \cdots x_{n+1}$
in $\Char(\Gamma_\w) \otimes \bQ \cong \bQ$
is greater than the degree of $\w$,
which is turn is greater
than the degree of any element
of the Jacobi ring
\begin{align}
 \Jac_{\w} \coloneqq \bC[x_1,\ldots,x_{n+1}] / (\partial_{x_1} \w, \ldots, \partial_{x_{n+1}} \w)
\end{align}
of $\w$,
and the proof of
\cite[Section 12.6, Theorem]{MR777682}
shows that
the polynomial
$\w+x_0 \cdots x_{n+1}$
considered as an element of
$\bC[\![x_0]\!][\![x_1,\ldots,x_{n+1}]\!]$
(i.e., a formal one-parameter deformation
of a formal germ of $\w$)
is right equivalent to $\w$
by a formal coordinate change
(i.e., there exists
$
 \varphi
  \in \Aut_{\bC[\![x_0]\!]} \bC[\![x_0]\!][\![x_1,\ldots,x_{n+1}]\!]
$
such that $\varphi^*(\w+x_0 \cdots x_{n+1}) = \w$).
%$\w$ defines a simple singularity and
The proof moreover shows that
one can choose $\varphi$
to be $\Gamma_\w$-equivariant,
which implies that
for any $i \in \{ 1, \ldots, n+1 \}$,
the coefficient $a_{i,m_1,\ldots,m_{n+1}}(x_0)$
of the expansion
$
 \varphi^* (x_i)
  = \sum_{m_1, \ldots, m_{n+1}=0}^\infty
   a_{i,m_1,\ldots,m_{n+1}}(x_0) x_1^{m_1} \cdots x_{n+1}^{m_{n+1}}
$
is a monomial in $x_0$,
since the degree of $x_0$ in $\Char(\Gamma_\w) \otimes \bQ$ is negative.
In particular,
one has
$
 \varphi
  \in \Aut_{\bC[x_0]} \bC[x_0][\![x_1,\ldots,x_{n+1}]\!].
%  \subset \Aut_{\bC[\![x_0]\!]} \bC[\![x_0]\!][\![x_1,\ldots,x_{n+1}]\!]
$
%so that
%$\w+x_0 \cdots x_n$
%considered as an element of
%the completion
%$\bC[x_0] [\![ x_1, \ldots, x_n ]\!]$
%along the critical locus
%is equivalent to
%$\w$ by a $\Gamma_\w$-equivariant automorphism of
%$\bC[x_0] [\![ x_1, \ldots, x_n ]\!]$,
This means that
that the formal completion of
$(\bA^{n+2}, \w+x_0 \cdots x_{n+1})$
along $x_1=\cdots=x_{n+1}=0$
is isomorphic to that of
$(\bA^{n+2}, \w)$
as a pair of a formal scheme and a regular function on it,
so that
\begin{align} \label{eq:completion}
 \mf \lb \bA^{n+2}, \Gamma_\w, \w+x_0 \cdots x_{n+1} \rb
  \simeq \mf \lb \bA^{n+2}, \Gamma_\w, \w \rb
\end{align}
by \cite[Theorem 2.10]{MR2735755},
and the quasi-equivalence \pref{eq:removing xyz} is proved.

\section{Hochschild cohomology of the derived $n$-preprojective algebra}
 \label{sc:HH}

\subsection{Hochschild cohomology via matrix factorizations}

We use the same notation as in \pref{sc:CY-completion},
and set
\begin{align}
 V \coloneqq \bC x_0 \oplus \bC x_1 \oplus \cdots \oplus \bC x_{n+1}.
\end{align}
For $\gamma \in \Gamma$,
let $V_\gamma$ be the subspace of $\gamma$-invariant elements in $V$,
$S_\gamma$ be the symmetric algebra of $V_\gamma$,
$\w_\gamma$ be the restriction of $\w$ to $\Spec S_\gamma$,
and
$N_\gamma$ be the $\Gamma$-stable complement of $V_\gamma$ in $V$
so that
%equipped with a $\Gamma$-action satisfying
$V \cong V_\gamma \oplus N_\gamma$
as a $\Gamma$-module.
Then
\cite{MR2824483,MR3084707,MR3108698,MR3270588}
(cf.~also \cite[Theorem 3.1]{1806.04345})
shows that
$
 \HH^t \lb \mf \lb \bA^{n+2}, \Gamma, \w \rb \rb
$
is isomorphic to
\begin{multline} \label{eq:HHmf}
 \bigoplus_{\substack{\gamma \in \ker \chi, \ l \geq 0 \\ t - \dim N_\gamma = 2u }}
  \lb
  H^{-2l}(d \w_\gamma) \otimes \Lambda^{\dim N_\gamma} N_\gamma^\dual \rb_{(u+\ell)\chi}
  \oplus \bigoplus_{\substack{\gamma \in \ker \chi, \ l \geq 0 \\ t - \dim N_\gamma = 2u+1}}
  \lb
  H^{-2l-1}(d \w_\gamma) \otimes \Lambda^{\dim N_\gamma} N_\gamma^\vee \rb_{(u+\ell+1) \chi}.
\end{multline}
Here $H^i(d \w_\gamma)$ is the $i$-th cohomology of the Koszul complex
\begin{align} \label{eq:Koszul}
 C^*(d \w_\gamma) \coloneqq \lc
 \cdots
 \to \Lambda^2 V_\gamma^\dual \otimes S_\gamma (-2 \chi)
  \to V_\gamma^\dual \otimes S_\gamma (-\chi)
  \to S_\gamma \rc,
\end{align}
%\begin{align}
% C(d\w) \coloneqq \lc \bigwedge^n \lb V^\dual(-\chi) \rb \otimes \cO_{\bA^{n+1}}
%  \to \cdots \to V^\dual(-\chi) \otimes \cO_{\bA^{n+1}}
%  \to \cO_{\bA^{n+1}} \rc,
%\end{align}
where the rightmost term $S_\gamma$ sits in cohomological degree 0, and
the differential is the contraction with
\begin{align}
 d \w_\gamma \in \lb V_\gamma \otimes S_\gamma \rb_{\chi}.
\end{align}

If $\w_\gamma$ has an isolated critical point at the origin,
then the cohomology of \pref{eq:Koszul}
is concentrated in degree 0,
so that only the summand
\begin{align} \label{eq:l=0}
 \lb \Jac_{\w_\gamma} \otimes \Lambda^{\dim N_\gamma} N_\gamma^\dual \rb_{u \chi}
\end{align}
with $l = 0$ in \pref{eq:HHmf} contributes to $\HH^{2 u + \dim N_\gamma}$.

If $V_\gamma$ contains $\bC x_0$,
then the Koszul complex
$C^*(d \w_\gamma)$
is isomorphic to the tensor product
of $C^*(d \w_\gamma')$
and the complex
$
 \lc \bC x_0^\dual \otimes \chi^\dual \otimes \bC[x_0] \to \bC[x_0] \rc
$
concentrated in cohomological degree $[-1,0]$
with the zero differential,
where $\w_\gamma'$ is the restriction of $\w$
to the complement $V_\gamma'$
of $\bC x_0$ in $V_\gamma$.
If $\w_\gamma'$ has an isolated critical point at the origin,
then $C^*(d \w_\gamma')$ is quasi-isomorphic to $\Jac_{\w_\gamma'}$
concentrated in cohomological degree $0$,
so that
only the summands
\begin{align} \label{eq:HH0_0}
 \lb \Jac_{\w_\gamma'}
  \otimes \bC[x_0]
  \otimes \Lambda^{\dim N_\gamma} N_\gamma^\dual \rb_{u \chi}
\end{align}
and
\begin{align} \label{eq:HH0_1}
 \lb \bC x_0^\dual \otimes \Jac_{\w_\gamma'}
  \otimes \bC[x_0]
  \otimes \Lambda^{\dim N_\gamma} N_\gamma^\dual \rb_{u \chi}
\end{align}
with $l=0$
in \pref{eq:HHmf}
contribute to $\HH^{2u+\dim N_\gamma}$ and
$\HH^{2u+\dim N_\gamma+1}$ respectively.

\begin{remark}
  Although \pref{eq:HHmf} may not look identical to \cite[Theorem 1.2]{MR3270588},
  the proof in \cite[Section 5]{MR3270588} actually shows \pref{eq:HHmf}.
  One way to think about \pref{eq:HHmf} is the following:
  If we set $\w = 0$ and forget $\Gamma$,
  then the Hochschild--Kostant--Rosenberg theorem gives
  a quasi-isomorphism
  of the Hochschild cochain complex of $\coh \bA^{n+2}$
  and
  \begin{align} \label{eq:HKR for Spec S}
    S \to V^\dual \otimes S \to \Lambda^2 V^\dual \otimes S \to \cdots
  \end{align}
  as complexes of $\bC$-vector spaces
  (which lifts to a quasi-isomorphism of $L_\infty$-algebras
  by the Kontsevich formality).
  If we introduce the potential $\w$,
  then the complex \pref{eq:HKR for Spec S} acquires an additional differential
  $\Lambda^i V^\dual \otimes S \to \Lambda^{i-1} V^\dual \otimes S$
  defined as the contraction with $d \w \in V \otimes S$,
  which decreases the cohomological grading by one
  so that the cohomological grading is collapsed to $\bZ/2\bZ$.
  The introduction of $\Gamma$ lifts the grading to $\bZ$ again
  and produces `twisted sectors' from the orbifold HKR theorem,
  leading to \pref{eq:HHmf};
  recall the isomorphism $(\chi) \simeq [2]$ of endofunctors of $\mf(\bA^{n+2},\Gamma,\w)$
  and the orbifold HKR theorem
  \begin{align}
    \HH^*([X/G])
    =
    \lb
    \bigoplus_{g \in G}
    \bigoplus_{p+q=*}
    H^{p-\codim X^g} \lb X^g, \Lambda^q T_{X^g} \otimes \Lambda^{\codim X^g} N_{X^g/X} \rb
    \rb_G
  \end{align}
  for global quotients
  appearing, e.g., in \cite[Corollary 1.17]{MR4003476}.
\end{remark}

\begin{remark}

The Hochschild cohomology
of a graded algebra $B$
(with no differential)
has a bigrading such that
\begin{align}
 \HH^{r+s} \lb B \rb^s
  \coloneqq \Ext_{B^{\mathrm{op}} \otimes B}^r \lb B, B[s] \rb.
\end{align}

When $B$ is the trivial extension algebra $B^{n}$
of a finite-dimensional algebra $A$,
by introducing a $\Gm$-action on $B^{n}$
such that $A$ has weight $0$
and $A^\dual[-n]$ has weight
%$1$,
$n$,
the $s$-grading on $\HH^* \lb B^{n} \rb$
can be described
%in a derived Morita invariant way
as
%$n$ times
the weight of the induced $\Gm$-action.

For any positive integer $m$,
the underlying ungraded algebra
of the trivial extension algebras
$B^{mn}$
is isomorphic to
$B^{n}$,
and only the cohomological gradings are different;
that of the former is $m$ times that of the latter.
It follows that one has an isomorphism
\begin{align}
 \HH^{r+ms} \lb B^{mn} \rb^{ms}
  \cong \HH^{r+s} \lb B^{n} \rb^s
\end{align}
of vector spaces for any positive integer $m$
such that the parities of $n$ and $mn$ are the same
(note that the signs in the Hochschild complex
depend on the parity of the cohomological grading).

When $Q$ is a Dynkin quiver,
one can transport the $\Gm$-action on $B_Q^{n}$
to $\scrG_Q^{n}$ through the Koszul duality \pref{eq:Koszul_duality},
so that $g$ for $g \in Q_1$ has weight $0$,
$g^*$ for $g \in Q_1$ has weight
%$-1$,
$-n$,
and $h_v$ for $v \in Q_0$ has weight
%$-1$.
$-n$.
This makes the isomorphism \pref{eq:hh} $\Gm$-equivariant,
so that the $\Gm$-weights on both sides agree.

Since $\w$ does not depend on $x_0$,
the $\Gm$-action on $\bA^{n+2}$
such that the weight of $x_i$ is
%$-1$
$-n$
for $i=0$
and $0$ for $i \in \{ 1, \ldots, n+1 \}$
keeps $\w$ invariant.
This induces a $\Gm$-action on $\mf(\bA^{n+2}, \Gamma, \w)$,
and hence on $B_Q^{n}$,
whose weight is $0$ on $A_Q$
and
%$1$
$n$
on $A_Q^\dual[-n]$
just as in \cite{1806.04345}.
This allows us to compute the $s$-grading on
$\HH^* \lb B_Q^{n} \rb$
as the $\Gm$-weight on \pref{eq:HHmf}.
This $\Gm$-action is mirror to the one
introduced in \cite{MR2929070} and
studied further for type A Milnor fibers in \cite{MR3033519}.
\end{remark}

\subsection{Type $A_\ell$}

Consider the case
\begin{align}
 \w = x_1^{\ell+1} + x_2^2 + \cdots + x_{n+1}^2 \in \bC[x_0, x_1,\ldots,x_{n+1}]
\end{align}
with
\begin{align}
 \Gamma
  &= \Gamma_\w
  \coloneqq \lc \gamma = (t_0, t_1, \ldots, t_{n+1}) \in (\Gm)^{n+2} \relmid
   t_1^{\ell+1} = t_2^2 = \cdots = t_{n+1}^2 = t_0 t_1 \cdots t_{n+1} \rc,
\end{align}
so that
$
 \ker \chi \cong \bmu_{\ell+1} \times \lb \bmu_2 \rb^{n}
$
and $\Char(\Gamma)$ is generated by
$\chi$ and $\chi_i = \deg x_i$ for $i \in \{ 0, \ldots, n+1 \}$
with relations
\begin{align}
 \chi = (\ell+1) \chi_1 = 2 \chi_2 = \cdots = 2 \chi_{n+1}
 = \chi_0 + \cdots + \chi_{n+1}.
\end{align}

\subsubsection{}

For any $\gamma \in \ker \chi$,
one has
\begin{align}
 \Jac_{\w_\gamma} =
\begin{cases}
 \bC[x_0] \otimes \bC[x_1]/(x_1^\ell) & \bC x_0 \oplus \bC x_1 \subset V_\gamma, \\
 \bC[x_0] & \bC x_0 \subset V_\gamma \text{ and } \bC x_1 \not \subset V_\gamma, \\
 \bC[x_1]/(x_1^\ell) & \bC x_0 \not \subset V_\gamma \text{ and } \bC x_1 \subset V_\gamma, \\
 \bC & \text{otherwise}.
\end{cases}
\end{align}
If we write an element of
$
 \Jac_{\w_\gamma} \otimes \Lambda^{\dim N_\gamma} N_\gamma^\dual
$
as
\begin{align}
 x_0^{k_0} x_1^{k_1}
  \otimes x_{j_1}^\dual \wedge x_{j_2}^\dual \wedge \ldots \wedge x_{j_s}^\dual,
\end{align}
where $k_0 = 0$ if $\bC x_0 \not \subset V_\gamma$
and $k_1 = 0$ if $\bC x_1 \not \subset V_\gamma$,
then its degree is given by
\begin{align}
 k_0 \chi_0 + k_1 \chi_1 - \chi_{j_1} - \cdots - \chi_{j_s},
\end{align}
which can be proportional to $\chi$
only if $V_\gamma$ is either $V$, $\bC x_0 \oplus \bC x_1$, $\bC x_0$, or $0$.
We now deal with each of these cases in turn.

\subsubsection{}

One has
$
 V_\gamma = V
$
if and only if $\gamma$ is the identity element.
The degree of
$
 x_0^{k_0} x_1^{k_1} \in \Jac_\w
$
is
\begin{align}
% k_0 (\chi - \chi_1 - \cdots - \chi_n) + k_1 \chi_1 =
  k_0 \chi - (k_0 - k_1) \chi_1 - k_0 \chi_2 - \cdots - k_0 \chi_{n+1},
\end{align}
which is proportional to $\chi$
if and only if $k_0$ is even and $\ell+1$ divides $k_0 - k_1$.
Such an element can be written as
\begin{align}
 \bsa_{k,m} \coloneqq x_0^{k+m(\ell+1)} x_1^k,
\end{align}
where
$
 k \in \{ 0, \ldots, \ell-1 \}
$
and
$
 m \in \bN
$
satisfies
\begin{itemize}
 \item
if $\ell$ is even, then the parities of $k$ and $m$ agree, and
 \item
if $\ell$ is odd, then $k$ is even.
\end{itemize}
Since
\begin{align}
 \deg \lb x_0^{k+m(\ell+1)} x_1^k \rb
  &= (k+m(\ell+1)) \chi - m \chi - \frac{1}{2} (k+m(\ell+1))n \chi \\
  &= \lb (k+m \ell) - \frac{1}{2} (k+m(\ell+1))n \rb \chi,
\end{align}
the element $x_0^{k+m(\ell+1)} x_1^k$
for such $(k,m)$ contributes
%$\bC(k+m(\ell+1))$
$\bC((k+m(\ell+1))n)$
to $\HH^t$ for $t = 2(k+m \ell) - (k+m(\ell+1))n$
by \pref{eq:HH0_0}.
Similarly,
for each such $(k,m)$,
the element
\begin{align}
 \bsalpha_{k,m} \coloneqq x_0^\dual \otimes x_0^{k+m(\ell+1)+1} x_1^k
   \in \bC x_0^\dual \otimes \Jac_{\w}
%  \in \lb \bC x_0^\dual \otimes \Jac_{\w} \otimes \chi^{\otimes \lb (k+m \ell) - \frac{1}{2} (k+m(\ell+1))n \rb} \rb^\Gamma,
\end{align}
contributes
$\bC((k+m(\ell+1))n)$ to $\HH^{t+1}$ for $t = 2(k+m \ell) - (k+m(\ell+1))n$
by \pref{eq:HH0_1}.

\subsubsection{}

One has
$
 V_\gamma = \bC x_0 \oplus \bC x_1
$
if and and only if
$n$ is even and
$
 \gamma = (1, 1, -1, \ldots, -1).
$
The degree of
\begin{align}
 x_0^{k_0} x_1^{k_1} \otimes x_2^\dual \wedge \cdots \wedge x_{n+1}^\dual
  \in  \Jac_{\w_\gamma} \otimes \Lambda^{\dim N_\gamma} N_\gamma^\dual
\end{align}
is given by
\begin{align}
 k_0 \chi + (k_1 - k_0) \chi_1 - (k_0+1) \chi_2 - \cdots - (k_0+1) \chi_{n+1},
\end{align}
which is proportional to $\chi$
if and only if $k_0$ is odd and $\ell+1$ divides $k_1 - k_0$.
Such an element can be written as
\begin{align}
 \bsa_{k,m} \coloneqq x_0^{k+m(\ell+1)} x_1^k \otimes x_2^\dual \wedge \cdots \wedge x_{n+1}^\dual,
\end{align}
where
$
 k \in \{ 0, \ldots, \ell-1 \}
$
and
$
 m \in \bN
$
satisfies
\begin{itemize}
 \item
if $\ell$ is even, then the parities of $k$ and $m$ differ, and
 \item
if $\ell$ is odd, then $k$ is odd.
\end{itemize}
Since the degree of this element is
\begin{align}
 \lb (k+m \ell) - \frac{1}{2} (k+m(\ell+1)+1)n \rb \chi,
\end{align}
each such $(k,m)$ contributes
%$\bC(k+m(\ell+1))$
$\bC((k+m(\ell+1))n)$
to $\HH^t$ for
\begin{align}
 t
  &= 2 \lb (k+m \ell) - \frac{1}{2}(k+m(\ell+1)+1)n \rb + \dim N_\gamma \\
  &= 2 (k+m \ell) - (k+m(\ell+1))n.
\end{align}
Similarly,
for each such $(k,m)$,
there is an element $\bsalpha_{k,m}$ contributing
%$\bC(k+m(\ell+1))$
$\bC((k+m(\ell+1))n)$
to $\HH^{t+1}$ for $t = 2(k+m \ell) - (k+m(\ell+1))n$.

\subsubsection{}

One has
$
 V_\gamma = \bC x_0
$
if and only if
both $\ell$ and $n$ are odd and
$
 \gamma = (1, -1, \ldots, -1).
$
The degree of
\begin{align}
 x_0^{k_0} \otimes x_1^\dual \wedge \cdots \wedge x_{n+1}^\dual
  \in  \Jac_{\w_\gamma} \otimes \Lambda^{\dim N_\gamma} N_\gamma^\dual
\end{align}
is given by
\begin{align}
 k_0 \chi - (k_0+1) \chi_1 - (k_0+1) \chi_2 - \cdots - (k_0+1) \chi_{n+1},
\end{align}
which is proportional to $\chi$
if and only if
%$k_0$ is odd and 
$\ell+1$ divides $k_0+1$.
Such an element can be written as
\begin{align}
 \bsb_m \coloneqq x_0^{m(\ell+1)-1} \otimes x_1^\dual \wedge \cdots \wedge x_{n+1}^\dual
\end{align}
for
$
 m \in \bN \setminus \{ 0 \}.
$
Since the degree of this element is
\begin{align}
 \lb (m \ell-1) - \frac{1}{2} m(\ell+1)n \rb \chi,
\end{align}
each such element contributes
$\bC((m(\ell+1)-1)n)$
to $\HH^t$ for
\begin{align}
 t
  &= 2 \lb (m \ell-1) - \frac{1}{2}m(\ell+1)n \rb + \dim N_\gamma \\
  &= (2 m \ell-1) - (m(\ell+1)-1)n.
\end{align}
Similarly,
for each $m \in \bN$,
the element
\begin{align}
 \bsbeta_m \coloneqq x_0^\dual \otimes x_0^{m(\ell+1)} \otimes x_1^\dual \wedge \cdots \wedge x_{n+1}^\dual
  \in  \bC x_0^\dual \otimes \Jac_{\w_\gamma} \otimes \Lambda^{\dim N_\gamma} N_\gamma^\dual
\end{align}
contributes
$\bC(m(\ell+1)-1)n)$
to $\HH^{t+1}$ for
$
 t = (2 m \ell-1) - (m(\ell+1)-1)n.
$

\subsubsection{}

For
$
 \gamma = (t_0, \ldots, t_{n+1}) \in \ker \chi,
$
one has $V_\gamma = 0$
if and only if $t_i \ne 1$ for all $i \in \{ 0, \ldots, n+1 \}$.
This is the case if and only if
$
 t_2 = \cdots = t_{n+1} = -1,
$
$
 t_1 \in \bsmu_{\ell+1} \setminus \{ 1 \},
$
and
\begin{align} \label{eq:t0_A}
 t_0 = (-1)^{n} t_1^{-1} \ne 1.
\end{align}
%Assume $t_1^{\ell+1} = 1$ and $t_1 \ne 1$.
If $n$ is odd, then \pref{eq:t0_A} holds if and only if $t_1 \ne -1$,
so that the number of such $\gamma$ is $\ell$ if $\ell$ is even,
and $\ell-1$ if $\ell$ is odd.
If $n$ is even, then \pref{eq:t0_A} always holds,
and the number of such $\gamma$ is $\ell$.
Each such $\gamma$ contributes
%$\bC(-1)$
$\bC(-n)$
to $\HH^{n}$.

\subsubsection{}

To sum up,
the Hochschild cohomology group has a basis
consisting of the following elements:
\begin{itemize}
 \item
$\bsa_{k,m}$ of degree
$
 2(k+m \ell) - (k+m (\ell+1))n
$
and weight
$
 -(k+m(\ell+1))n,
$
and
$\bsalpha_{k,m}$ of degree
$
 2(k+m \ell) - (k+m (\ell+1))n+1
$
and weight
%$-(k+m(\ell+1))$,
$
 -(k+m(\ell+1))n,
$
where
\begin{itemize}
 \item
 if $n$ is even, then
 $(k,m)$ runs over $\{ 0, \ldots, \ell-1 \} \times \bN$,
 and
 \item
if $n$ is odd, then $(k,m)$ runs over those pairs in $\{ 0, \ldots, \ell-1 \} \times \bN$
for which
\begin{itemize}
 \item
the parities of $k$ and $m$ agree, if $\ell$ is even,
 \item
$k$ is even, if $\ell$ is odd,
\end{itemize}
\end{itemize}
 \item
if both $n$ and $\ell$ are odd, then
%in addition to the above, one has
\begin{itemize}
 \item
$\bsb_m$ of degree
$
 2m \ell-1-(m(\ell+1)-1)n
$
and weight
$
 -(m(\ell+1)-1)n
$
for $m \in \bN \setminus \{ 0 \}$,
and
 \item
$\bsbeta_{m}$ of degree
$
 2 m \ell-(m(\ell+1)-1)n
$
and weight
$
 -(m(\ell+1)-1)n
$
for $m \in \bN$, and
\end{itemize}
 \item
$\bss_h$ of degree $n$ and weight
%$1$
$n$,
where $h$ runs over
\begin{itemize}
 \item
$\{ 1, 2, \ldots, \ell-1 \}$ if both $\ell$ and $n$ are odd, and
 \item
$\{ 1, 2, \ldots, \ell \}$ otherwise.
\end{itemize}
\end{itemize}

\subsubsection{}

As an example,
consider the case $\ell=1$. Note that the $A_1$-Milnor fiber is symplectomorphic to the cotangent bundle $T^*S^n$. 
The Hochschild cohomology group
in this case
is spanned by
\begin{itemize}
 \item
$\bsa_{0,m}$ for $m \in \bN$ of
degree $-2m(n-1)$ and weight $-2mn$,
 \item
$\bsalpha_{0,m}$ for $m \in \bN$ of
degree $-2m(n-1)+1$ and weight $-2mn$,
\end{itemize}
and, if $n$ is odd, in addition to the above,
\begin{itemize}
 \item
%if $n$ is even, then
$\bsb_m$ for $m \in \bN \setminus \{ 0 \}$ of
degree $-(2m-1)(n-1)$ and weight $-(2m-1)n$,
 \item
%if $n$ is even, then
$\bsbeta_m$ for $m \in \bN$ of
degree $-(2m-1)(n-1)+1$ and weight $-(2m-1)n$,
\end{itemize}
and, if $n$ is even, in addition to the above,
\begin{itemize}
 \item
$\bss_1$ of degree $n$ and weight $n$.
\end{itemize}
This is consistent with the isomorphism
\begin{align}
 \SH^*(T^* S^{n})
  \cong H_{n-*}(\scrL S^{n}),
\end{align}
which is a special case of the isomorphism
between the symplectic cohomology of the cotangent bundle
and the homology of the free loop space
\cite[Theorem 3.1]{1805.01316}
(see e.g.~\cite[Theorem 2]{MR2039760}
for the homology of the free loop space of spheres).

Another example is the case when $n=2$ and $\ell$ is arbitrary.
In this case, $\SH^*(\wv^{-1}(1))$ was computed in \cite{MR3692968}
as a bigraded ring.
This is compatible with the computation given here.

\subsection{Type $D_\ell$}
 \label{sc:D}

The Berglund--H\"ubsch transform of the invertible polynomial
\begin{align}
  \check{\bfv} = y_1^{\ell-1} + y_1 y_2^2 + y_3^2 + \cdots + y_{n+1}^2
\end{align}
defining the $D_\ell$-singularity is given by
\begin{align}
  \bfv = y_1^{\ell-1} y_2 + y_2^2 + \cdots + y_{n+1}^2,
\end{align}
and one has
\begin{align}
  \Gamma_\bfv = \lc \gamma = (t_1, \ldots, t_{n+1}) \in (\Gm)^{n+1} \relmid
  t_1^{\ell-1} t_2 = t_2^2 = \cdots = t_{n+1}^2 \rc.
\end{align}
By completing the square and rescaling,
one has
\begin{align}
  \bfv(y) = \w(x(y))
\end{align}
where
\begin{align} \label{eq:coordinate transform for type D}
  x_1 = (-1/4)^{1/(2n-2)} y_1, \ 
  x_2 = y_2 + \frac{1}{2} y_1^{\ell-1}, \ 
  x_3 = y_3, \ \ldots, \ 
  x_{n+1} = y_{n+1}
\end{align}
% $
% x_1 = (-1/4)^{1/(2n-2)} y_1,
% $
% $
% x_2 = y_2 + \frac{1}{2} y_1^{\ell-1},
% $
% $
% x_i = y_i
% $
% for
% $
% i = 3, \ldots, n+1,
% $
and
\begin{align}
  \w = x_1^{2 \ell - 2} + x_2^2 + \cdots + x_{n+1}^2.
\end{align}
Although 
the change of variables \pref{eq:coordinate transform for type D}
is neither linear nor diagonal,
the induced action of $\Gamma_{\bfv}$
on $\Spec \bC[x_1,\ldots,x_{n+1}]$
remains linear and diagonal,
so that one can identify $\Gamma_{\bfv}$
with a proper subgroup of $\Gamma_{\w}$.

Therefore,
we will work with
\begin{align}
 \w = x_1^{2\ell-2} + x_2^2 + \cdots + x_{n+1}^2 \in \bC[x_0, x_1,\ldots,x_{n+1}]
\end{align}
with the non-maximal group
\begin{align}
 \Gamma
  = \lc \gamma = (t_0, t_1, \ldots, t_{n+1}) \in (\Gm)^{n+2} \relmid
   t_1^{\ell-1} t_2 = t_2^2 = \cdots = t_{n+1}^2 = t_0 t_1 \cdots t_{n+1} \rc.
\end{align}
One has
$
 \ker \chi \cong \bmu_{2\ell-2} \times \lb \bmu_2 \rb^{n-1}
$
and $\Char(\Gamma)$ is generated by
$\chi$ and $\chi_i = \deg x_i$ for $i \in \{ 0, \ldots, n+1 \}$
with relations
\begin{align} \label{eq:D_rel}
 \chi = (\ell-1) \chi_1 + \chi_2 = 2 \chi_2 = \cdots = 2 \chi_{n+1} = \chi_0 + \cdots + \chi_{n+1}.
\end{align}
% By a change of coordinates,
% this is equivalent to
% \begin{align}
%  \w' = x_1^{\ell-1} x_2 + x_2^2 + \cdots + x_{n+1}^2
% \end{align}
% with
% $
%  \Gamma = \Gamma_{\w'},
% $
% whose Berglund--H\"ubsch transpose
% \begin{align}
%  \wv' = x_1^{\ell-1} + x_1 x_2^2 + x_3^2 + \cdots + x_{n+1}^2
% \end{align}
% defines the $D_\ell$-singularity.
The relations \pref{eq:D_rel} imply
\begin{align}
 \chi_2 &= \chi - (\ell-1) \chi_1, \\
 \chi &= (2\ell-2) \chi_1, \\
 \chi_0 &= \chi - \chi_1 - \cdots - \chi_n \\
 &= (\ell-2) \chi_1 - \chi_3 - \cdots - \chi_{n+1}.
\end{align}

\subsubsection{}

For any
$
\gamma
= (t_0,\ldots,t_{n+1})
\in \ker \chi,
$
one has
\begin{align}
 \Jac_{\w_\gamma} =
\begin{cases}
 \bC[x_0] \otimes \bC[x_1]/(x_1^{2\ell-3}) & \bC x_0 \oplus \bC x_1 \subset V_\gamma, \\
 \bC[x_0] & \bC x_0 \subset V_\gamma \text{ and } \bC x_1 \not \subset V_\gamma, \\
 \bC[x_1]/(x_1^{2\ell-3}) & \bC x_0 \not \subset V_\gamma \text{ and } \bC x_1 \subset V_\gamma, \\
 \bC & \text{otherwise}.
\end{cases}
\end{align}
If we write an element of
$
 \Jac_{\w_\gamma} \otimes \Lambda^{\dim N_\gamma} N_\gamma^\dual
$
as
\begin{align}
 x_0^{k_0} x_1^{k_1}
  \otimes x_{j_1}^\dual \wedge x_{j_2}^\dual \wedge \ldots \wedge x_{j_s}^\dual,
\end{align}
where $k_0 = 0$ if $\bC x_0 \not \subset V_\gamma$
and $k_1 = 0$ if $\bC x_1 \not \subset V_\gamma$,
then its degree is given by
\begin{align}
 k_0 \chi_0 + k_1 \chi_1 - \chi_{j_1} - \cdots - \chi_{j_s},
\end{align}
which can be proportional to $\chi$
only if
\begin{align}
V_\gamma \cap \lb \bC x_3 \oplus \cdots \oplus \bC x_{n+1} \rb
\text{ is either }
\bC x_3 \oplus \cdots \oplus \bC x_{n+1}
\text{ or }
0,
\end{align}
that is,
\begin{align} \label{eq:condition for type D}
 t_3=\cdots=t_{n+1} = \pm 1.
\end{align}
% \begin{itemize}
%  \item
% $
%  V_\gamma \cap \lb \bC x_3 \oplus \cdots \oplus \bC x_{n+1} \rb
% $
% is either
% $
%  \bC x_3 \oplus \cdots \oplus \bC x_{n+1}
% $
% or $0$.
% \end{itemize}
We will assume this condition
for the rest of \pref{sc:D}.

\subsubsection{}

One has
$
\gamma = (t_0,\ldots,t_{n+1}) \in \ker \chi
$
if and only if
\begin{align}
t_1^{\ell-1} t_2 = t_2^2 = \cdots = t_{n+1}^2 = t_0 t_1 \cdots t_{n+1} = 1.
\end{align}

If $t_0 = 1$,
then one has $t_2^2 = \cdots = t_{n+1}^2 = t_1 \cdots t_{n+1} = 1$,
so that
$t_1 = (t_2 \cdots t_{n+1})^{-1} = \pm 1$.

If $t_1 = 1$, then one has
$
t_2 = t_2^2,
$
so that
$t_2 = 1$
and
$
t_3^2 = \cdots = t_{n+1}^2 = t_0 t_3 \cdots t_{n+1} = 1.
$
Under the assumption \pref{eq:condition for type D},
one has $t_0 = 1$ if and only if $(t_3)^{n-1} = 1$,
that is, $t_3 = 1$ or $n$ is even.

% \item
If $t_1=-1$, then $t_2 = (-1)^{\ell-1}$,
and one has $t_0 = 1$ if and only if $(-1)^{\ell} t_3^{n-1} = 1$.
%that is, the parities of $\ell$ and $n$ are the same.
%\end{itemize}
It follows that
\begin{itemize}
 \item
$V_\gamma$ contains $\bC x_0$
if and only if
\begin{itemize}
 \item
$\gamma = (1,\ldots,1)$, where $V_\gamma = V$,
 \item
$\gamma = (1,1,1,-1,\ldots,-1)$ with odd $n$,
where $V_\gamma = \bC x_0 \oplus \bC x_1 \oplus \bC x_2$,
 \item
$\gamma = (1,-1,-1,1,\ldots,1)$ with even $\ell$,
where $V_\gamma = \bC x_0 \oplus \bC x_3 \oplus \cdots \oplus \bC x_{n+1}$,
 \item
$\gamma = (1,-1,-1,-1,\ldots,-1)$ with even $\ell$ and odd $n$,
where $V_\gamma = \bC x_0$,
 \item
$\gamma = (1,-1,1,-1,\ldots,-1)$ with odd $\ell$ and even $n$,
where $V_\gamma = \bC x_0 \oplus \bC x_2$.
\end{itemize}
\end{itemize}

\subsubsection{}

Note for later use that
the smallest positive integer $k$ such that the degree of $x_0^k$
is proportional to $\chi$ is $2 \ell-2$.
One has
\begin{align}
 \deg x_0^{2\ell-2}
  &= (2\ell-2)(\chi-\chi_1-\cdots-\chi_{n+1}) \\
%  &= (2\ell-2)\chi-\chi-(\ell-1)n \chi \\
  &= \lb (2\ell-3)-(\ell-1) n \rb \chi.
\end{align}

\subsubsection{}

One has
$
 V_\gamma = V
$
if and only if $\gamma$ is the identity element.
The degree of
$
 x_0^{k_0} x_1^{k_1} \in \Jac_\w
$
is
\begin{align}
 k_0 \chi - (k_0-k_1) \chi_1 - k_0 \chi_2 - \cdots - k_0 \chi_{n+1},
\end{align}
which is proportional to $\chi$
if and only if $k_0$ is even and $2\ell-2$ divides $k_0 - k_1$.
Such an element can be written as
\begin{align}
 \bsa_{k,m} \coloneqq x_0^{2k+(2\ell-2)m} x_1^{2k}
\end{align}
for
$
 (k,m) \in \{ 0, \ldots, \ell-2 \} \times \bN
$
which contributes
$
 \bC((2k+(2\ell-2)m)n)
$
to $\HH^t$ for
$
 t = 4k+(4\ell-6)m - (2k+(2\ell-2)m)n
$
since
\begin{align}
 \deg x_0^{2k} x_1^{2k}
  &= (2k - kn) \chi.
\end{align}
Similarly,
for each $(k,m) \in \{ 0, \ldots, \ell-2 \} \times \bN$,
there is an element $\bsalpha_{k,m}$ contributing
$
 \bC((2k+(2\ell-2)m)n)
$
to $\HH^t$ for
$
 t = 4k+1+(4\ell-6)m - (2k+(2\ell-2)m)n.
$

\subsubsection{}

One has
$
 V_\gamma = \bC x_0 \oplus \bC x_1 \oplus \bC x_2
$
if and only if $\gamma = (1,1,1,-1,\ldots,-1)$ and $n$ is odd.
The degree of
$
 x_0^{k_0} x_1^{k_1} \otimes x_3^\dual \wedge \cdots \wedge x_{n+1}^\dual
  \in \Jac_{\w_\gamma} \otimes \Lambda^{\dim N_\gamma} N_\gamma^\dual
$
is
\begin{align}
 k_0 \chi - (k_0 - k_1) \chi_1 - k_0 \chi_2 - (k_0+1) \chi_3 - \cdots - (k_0+1) \chi_{n+1},
\end{align}
which is proportional to $\chi$
if and only if $k_0$ is odd and $2\ell-2$ divides $k_0 - k_1-(\ell-1)$.
%since
%\begin{align}
% \chi_2 \equiv (\ell-1) \chi_1 \mod \chi.
%\end{align}
Such an element can be written as
\begin{align}
 \bsb_{k,m} \coloneqq
  x_0^{k+\ell-1+(2\ell-2)m} x_1^{k} \otimes x_3^\dual \wedge \cdots \wedge x_{n+1}^\dual
\end{align}
for
\begin{align}
 (k,m) \in \lc (k,m) \in \{ 0, \ldots, 2 \ell-4 \} \times \bZ \relmid
  k+\ell \text{ is even and } k+\ell-1+m(2\ell-2) \ge 0 \rc.
\end{align}
It contributes
$
 \bC((k+\ell-1+(2\ell-2)m)n)
$
to $\HH^t$ for
\begin{align}
 t
  &= 2 \deg \lb x_0^{k+\ell-1+(2\ell-2)m} x_1^k \otimes x_3^\dual \wedge \cdots \wedge x_{n+1}^\dual \rb / \chi + \dim N_\gamma \\
 &= 2k+2\ell-3+(4\ell-6)m-(k+\ell-1-(2\ell-2)m)n,
\end{align}
since
\begin{align}
 &\deg \lb x_0^{k+\ell-1} x_1^k \otimes x_3^\dual \wedge \cdots \wedge x_{n+1}^\dual \rb \\
 &\qquad= (k+\ell-1) \chi - (\ell-1) \chi_1 - (k+\ell-1) \chi_2 - (k+\ell) \chi_3 - \cdots - (k+\ell) \chi_{n+1} \\
% &\qquad= (k+\ell-1) \chi - ((\ell-1) \chi_1 + \chi_2) + 2 \chi_2 - (k+\ell) \chi_2 - \cdots - (k+\ell) \chi_{n+1} \\
% &\qquad= (k+\ell-1) \chi - \chi + \chi - \frac{1}{2} (k+\ell) n \chi \\
 &\qquad = \lb k+\ell-1 -\frac{1}{2}(k+\ell)n \rb \chi.
\end{align}
Similarly,
for each
\begin{align}
 (k,m) \in \lc (k,m) \in \{ 0, \ldots, 2 \ell-4 \} \times \bZ \relmid
  k+\ell \text{ is even and } k+\ell+m(2\ell-2) \ge 0 \rc,
\end{align}
the element
\begin{align}
 \bsbeta_{k,m} \coloneqq
 x_0^\dual \otimes x_0^{k+\ell+(2\ell-2)m} x_1^{k} \otimes x_3^\dual \wedge \cdots \wedge x_{n+1}^\dual
   \in \lb (\bC x_0)^\dual \otimes \Jac_{\w_\gamma} \otimes \Lambda^{\dim N_\gamma} N_\gamma^\dual \rb^\Gamma
\end{align}
contributes
$
 \bC(k+\ell-1+(2\ell-2)m)
$
to $\HH^t$ for
\begin{align}
 t = 2k+2\ell-2+(4\ell-6)m-(k+\ell-1-(2\ell-2)m)n.
\end{align}

\subsubsection{}

One has
$
 V_\gamma = \bC x_0 \oplus \bC x_3 \oplus \cdots \oplus \bC x_{n+1}
$
if and only if $\gamma = (1,-1,-1,1,\ldots,1)$
with even $\ell$.
An element of
$
 \Jac_{\w_\gamma} \otimes \Lambda^{\dim N_\gamma} N_\gamma^\dual
$
whose degree is proportional to $\chi$
can be written as
\begin{align}
 \bsc_m \coloneqq x_0^{\ell-2 + (2\ell-2)m} \otimes x_1^\dual \wedge x_2^\dual
\end{align}
for
$
 m \in \bN,
$
which contributes
$\bC((\ell-2+(2\ell-2)m)n)$
to $\HH^t$ for
\begin{align}
 t
  &= 2 \deg \lb x_0^{\ell-2} \otimes x_1^\dual \wedge x_2^\dual \rb / \chi + \dim N_\gamma \\
  &= 2 \ell-4+(4\ell-6)m-(\ell-2+(2\ell-2)m)n
\end{align}
since
\begin{align}
 \deg \lb x_0^{\ell-2} \otimes x_1^\dual \wedge x_2^\dual \rb
%  &= (\ell-2)(\chi-\chi_1-\cdots-\chi_{n+1})-\chi_1-\chi_2 \\
%  &= (\ell-2) \chi - (\ell-1) \chi_1 - \chi_2 - \frac{1}{2} (\ell-2) n \chi \\
  &= \lb \ell-3 - \frac{1}{2} (\ell-2)n \rb \chi.
\end{align}
Similarly,
for each $m \in \bN$,
there is an element $\bsgamma_m$
contributing
$\bC((\ell-2+(2\ell-2)m)n)$
to $\HH^t$ for
\begin{align}
 t
  &= 2\ell-3+(4\ell-6)m-(\ell-2+(2\ell-2)m)n.
\end{align}

\subsubsection{}

One has
$
 V_\gamma = \bC x_0
$
if and only if
$\ell$ is even, $n$ is odd, and
$
 \gamma = (1, -1, \ldots, -1) \in \ker \chi.
$
The degree of
\begin{align}
 x_0^{k_0} \otimes x_1^\dual \wedge \cdots \wedge x_{n+1}^\dual
  \in  \Jac_{\w_\gamma} \otimes \Lambda^{\dim N_\gamma} N_\gamma^\dual
\end{align}
is given by
\begin{align}
 k_0 \chi - (k_0+1) \chi_1 - (k_0+1) \chi_2 - \cdots - (k_0+1) \chi_{n+1},
\end{align}
which is proportional to $\chi$
if and only if
%$k_0$ is odd and 
$2\ell-2$ divides $k_0+1$.
Such an element can be written as
\begin{align}
 \bsd_m \coloneqq
 x_0^{-1+m(2\ell-2)} \otimes x_1^\dual \wedge \cdots \wedge x_{n+1}^\dual
\end{align}
for
$
 m \in \bN \setminus \{ 0 \}.
$
Since
\begin{align}
 \deg \lb x_0^{-1} \otimes x_1^\dual \wedge \cdots \wedge x_{n+1}^\dual \rb
  = -\chi,
\end{align}
each such element contributes
$\bC((-1+(2\ell-2)m)n)$
to $\HH^t$ for
\begin{align}
 t
  &= 2 \deg \lb x_0^{-1+(2\ell-2)m} \otimes x_1^\dual \wedge \cdots \wedge x_{n+1}^\dual \rb / \chi
   + \dim N_\gamma \\
  &= -1+(4\ell-6)m - (-1+(2\ell-2)m)n.
\end{align}
Similarly,
for each $m \in \bN$,
the element
\begin{align}
 \bsdelta_m \coloneqq
 x_0^\dual \otimes x_0^{m(2\ell-2)} \otimes x_1^\dual \wedge \cdots \wedge x_{n+1}^\dual
  \in  \bC x_0^\dual \otimes \Jac_{\w_\gamma} \otimes \Lambda^{\dim N_\gamma} N_\gamma^\dual
\end{align}
contributes
$\bC((-1+(2\ell-2)m)n)$
to $\HH^t$ for
\begin{align}
 t
  &= (4\ell-6)m - (-1+(2\ell-2)m)n.
\end{align}

\subsubsection{}

One has
$
 V_\gamma = \bC x_0 \oplus \bC x_2
$
if and only if
$\ell$ is odd,
$n$ is even, and
$
 \gamma = (1, -1, 1, -1, \ldots, -1).
$
The degree of
\begin{align}
 x_0^{k_0} \otimes x_1^\dual \wedge x_3^\dual \wedge \cdots \wedge x_{n+1}^\dual
  \in  \Jac_{\w_\gamma} \otimes \Lambda^{\dim N_\gamma} N_\gamma^\dual
\end{align}
is given by
\begin{align}
 k_0 \chi_0 - \chi_1 - \chi_3 - \cdots - \chi_{n+1}
% &\qquad= k_0 (\ell-2) \chi_1 - \chi_1 - (k_0+1) \chi_3 -\cdots- (k_0+1) \chi_{n+1} \\
 &= (k_0 (\ell-2) - 1) \chi_1 - (k_0+1) \chi_3 -\cdots- (k_0+1) \chi_{n+1},
\end{align}
which is proportional to $\chi$
if and only if $k_0$ is odd and $2\ell-2$ divides $k_0(\ell-2)-1$.
Such an element can be written as
\begin{align}
 \bse_m \coloneqq x_0^{\ell-2+(2\ell-2)m}
  \otimes x_1^\dual \wedge x_3^\dual \wedge \cdots \wedge x_{n+1}^\dual
\end{align}
for
$
 m \in \bN,
$
which contributes
$
 \bC((\ell-2+(2\ell-2)m)n)
$
to $\HH^t$ for
\begin{align}
 t
%  &= 2\ell-4+(4\ell-6)m-(\ell-1+(2\ell-2)m)n + n \\
  &= 2\ell-4+(4\ell-6)m-(\ell-2+(2\ell-2)m)n
\end{align}
since
\begin{align}
 \deg \lb x_0^{\ell-2}
  \otimes x_1^\dual \wedge x_3^\dual \wedge \cdots \wedge x_{n+1}^\dual \rb
%  &\qquad= ((\ell-2)^2-1)\chi_1-\frac{1}{2}(\ell-1)(n-1)\chi \\
%  &\qquad= (\ell^2-4\ell+3)\chi_1-\frac{1}{2}(\ell-1)(n-1)\chi \\
%  &\qquad= (\ell-1)(\ell-3)\chi_1-\frac{1}{2}(\ell-1)(n-1)\chi \\
%  &\qquad= \frac{1}{2}(\ell-3)\chi-\frac{1}{2}(\ell-1)(n-1)\chi \\
  &= \frac{1}{2} \lb 2\ell-4-(\ell-1)n \rb \chi.
\end{align}
Similarly,
for each
$
 m \in \bN,
$
there is an element $\bsepsilon_m$ contributing
$
 \bC((\ell-2+(2\ell-2)m)n)
$
to $\HH^t$ for
\begin{align}
 t &= 2\ell-3+(4\ell-6)m-(\ell-2+(2\ell-2)m)n.
\end{align}

\subsubsection{}

Now we move on to the case when
$\bC x_0 \not \subset V_\gamma$.
We divide it into three cases:
\begin{itemize}
  \item $\bC x_1 \subset V_\gamma$.
  \item $\bC x_1 \not \subset V_\gamma$ and $V_\gamma \neq 0$.
  \item $V_\gamma=0$.
\end{itemize}

\subsubsection{}

Set
$
 \zeta \coloneqq \exp \lb 2 \pi \sqrt{-1}/(2\ell-2) \rb.
$
For a given $\gamma = (t_0, \ldots, t_{n+1}) \in \ker \chi$,
we write $t_1 = \zeta^p$ for $p \in \{ 0, \ldots, 2 \ell-3 \}$.
Then one has $t_2 = (-1)^p$,
so that
\begin{itemize}
 \item
$V_\gamma$ contains $\bC x_1$ if and only if $p=0$, and
 \item
$V_\gamma$ contains $\bC x_2$ if and only if $p$ is even.
\end{itemize}

\subsubsection{}

If $\bC x_0 \not \subset V_\gamma$
and $\bC x_1 \subset V_\gamma$,
then one has that
$\gamma = (-1,1,1,-1,\ldots,-1)$,
$n$ is even, and
$V_\gamma = \bC x_1 \oplus \bC x_2$.
The element
\begin{align}
 x_1^{\ell-2} \otimes x_0^\dual \wedge x_3^\dual \wedge \cdots \wedge x_{n+1}^\dual
\end{align}
has degree
\begin{align}
 (\ell-2) \chi_1 - \chi_0 - \chi_3 - \cdots - \chi_{n+1}
  &= 0,
\end{align}
so that it contributes $\bC(-n)$ to $\HH^t$
for
$
 t = \dim N_\gamma = n,
$
and this is the only contribution.

\subsubsection{}

If $\bC x_0 \not \subset V_\gamma$,
$\bC x_1 \not \subset V_\gamma$,
and
$V_\gamma \neq 0$,
then $V_\gamma$ is either
$V_\gamma = \bC x_2$,
$\bC x_2 \oplus \cdots \oplus \bC x_{n+1}$, or
$\bC x_3 \oplus \cdots \oplus \bC x_{n+1}$.
No such $\gamma$ does not contribute to $\HH^*$,
since
$
 \Jac_{\w_\gamma} \otimes \Lambda^{\dim N_\gamma} N_\gamma^\dual
$
is spanned by a single element,
whose degree is not proportional to $\chi$.

\subsubsection{}

One has $V_\gamma = 0$
%if and only if $t_i \ne 1$ for all $i \in \{ 0, \ldots, n \}$.
%This is the case
if and only if
$
 t_3 = \cdots = t_{n+1} = -1,
$
$
 t_1 = \zeta^{2m+1}
$
for $m \in \{ 0, \ldots, \ell-2 \}$,
and
\begin{align} \label{eq:t0_D}
 t_0 = (-1)^{n} \zeta^{-2m-1} \ne 1.
\end{align}
The number of such $\gamma$ is
$\ell-2$ if $\ell$ is even and $n$ is odd,
and $\ell-1$ otherwise.
Each such $\gamma$ contributes $\bC(-n)$ to $\HH^{n}$.

\subsubsection{}

To sum up,
the Hochschild cohomology group has a basis
consisting of the following elements:
\begin{itemize}
 \item
$\bsa_{k,m}$ of degree
$
 4k+(4\ell-6)m-(2k+(2\ell-2)m)n
$
and weight
$
 - (2k+(2\ell-2)m)n
$
for
$
 (k,m) \in \{ 0, \ldots, \ell-2 \} \times \bN,
$
 \item
$\bsalpha_{k,m}$ of degree
$
 4k+1+(4\ell-6)m-(2k+(2\ell-2)m)n
$
and weight
$
 - (2k+(2\ell-2)m)n
$
for
$
 (k,m) \in \{ 0, \ldots, \ell-2 \} \times \bN,
$
 \item
if $n$ is odd,
$\bsb_{k,m}$ of degree
$
 2k+2\ell-3+(4\ell-6)m-(k+\ell-1-(2\ell-2)m)n
$
and weight
$
 -(k+\ell-1+(2\ell-2)m)n
$
for
$
 \lc (k,m) \in \{ 0, \ldots, 2 \ell-4 \} \times \bZ \relmid
  k+\ell \text{ is even and } k+\ell-1+m(2\ell-2) \ge 0 \rc,
$
 \item
if $n$ is odd,
$\bsbeta_{k,m}$ of degree
$
 2k+2\ell-2+(4\ell-6)m-(k+\ell-1-(2\ell-2)m)n
$
and weight
$
 -(k+\ell-1-(2\ell-2)m)n
$
for
$
 \lc (k,m) \in \{ 0, \ldots, 2 \ell-4 \} \times \bZ \relmid
  k+\ell \text{ is even and } k+\ell+m(2\ell-2) \ge 0 \rc,
$
 \item
if $\ell$ is even,
$\bsc_m$ of degree
$
 2 \ell-4+(4\ell-6)m-(\ell-2+(2\ell-2)m)n
$
and weight
$
 -(\ell-2+(2\ell-2)m)n
$
for
$
 m \in \bN,
$
 \item
if $\ell$ is even,
$\bsgamma_m$ of degree
$
 2 \ell-3+(4\ell-6)m-(\ell-2+(2\ell-2)m)n
$
and weight
$
 -(\ell-2+(2\ell-2)m)n
$
for
$
 m \in \bN,
$
 \item
if $\ell$ is even and $n$ is odd,
$\bsd_m$ of degree
$
 -1+(4\ell-6)m - (-1+(2\ell-2)m)n
$
and weight
$
 - (-1+(2\ell-2)m)n
$
for $m \in \bN \setminus \{ 0 \}$,
 \item
if $\ell$ is even and $n$ is odd,
$\bsdelta_m$ of degree
$
 (4\ell-6)m - (-1+(2\ell-2)m)n
$
and weight
$
 - (-1+(2\ell-2)m)n
$
for $m \in \bN$,
 \item
if $\ell$ is odd and $n$ is even,
$\bse_m$ of degree
$
 2\ell-4+(4\ell-6)m-(\ell-2+(2\ell-2)m)n
$
and weight
$
 -(\ell-2+(2\ell-2)m)n
$
for $m \in \bN$,
 \item
if $\ell$ is odd and $n$ is even,
$\bsepsilon_m$ of degree
$
 2\ell-3+(4\ell-6)m-(\ell-2+(2\ell-2)m)n
$
and weight
$
 -(\ell-2+(2\ell-2)m)n
$
for $m \in \bN$,
and
 \item
$\bss_h$ of degree $n$ and weight $n$,
where $h$ runs over a set consisting of
\begin{itemize}
 \item
$\ell-2$ elements if $\ell$ is even and $n$ is odd,
 \item
$\ell-1$ elements if both $\ell$ and $n$ are odd, and
 \item
$\ell$ elements otherwise.
\end{itemize}
\end{itemize}

\subsection{Type $E_6$}
 \label{sc:E6}

Consider the case
\begin{align}
 \w = x_1^4 + x_2^3 + x_3^2 + \cdots + x_{n+1}^2 \in \bC[x_0, x_1,\ldots,x_{n+1}]
\end{align}
with
\begin{align}
 \Gamma
  = \Gamma_\w
  \coloneqq \lc \gamma = (t_0, t_1, \ldots, t_{n+1}) \in (\Gm)^{n+2} \relmid
   t_1^4 = t_2^3 = t_3^2 = \cdots = t_{n+1}^2 = t_0 t_1 \cdots t_{n+1} \rc,
\end{align}
so that
$
 \ker \chi \cong \bmu_4 \times \bmu_3 \times \lb \bmu_2 \rb^{n-1}
$
and $\Char(\Gamma)$ is generated by
$\chi$ and $\chi_i = \deg x_i$ for $i \in \{ 0, \ldots, n+1 \}$
with relations
\begin{align}
 \chi = 4 \chi_1 = 3 \chi_2 = 2 \chi_3 = \cdots = 2 \chi_{n+1} = \chi_0 + \cdots + \chi_{n+1}.
\end{align}

\subsubsection{}

For any $\gamma \in \ker \chi$,
one has
\begin{align}
 \Jac_{\w_\gamma} \cong
\begin{cases}
 \bC[x_0] & \bC x_0 \subset V_\gamma \\
 \bC & \bC x_0 \not \subset V_\gamma
\end{cases}
 \otimes
\begin{cases}
 \bC[x_1]/(x_1^3) & \bC x_1 \subset V_\gamma \\
 \bC & \bC x_1 \not \subset V_\gamma
\end{cases}
 \otimes
\begin{cases}
 \bC[x_2]/(x_2^2) & \bC x_2 \subset V_\gamma \\
 \bC & \bC x_2 \not \subset V_\gamma.
\end{cases}
\end{align}
%the Jacobi ring $\Jac_{\w_\gamma}$ is isomorphic
%to the tensor product of
%\begin{align}
%\begin{cases}
% \bC[x_0] & \bC x_0 \subset V_\gamma, \\
% \bC & \bC x_0 \not \subset V_\gamma,
%\end{cases}
%\end{align}
%\begin{align}
%\begin{cases}
% \bC[x_1]/(x_1^3) & \bC x_1 \subset V_\gamma, \\
% \bC & \bC x_1 \not \subset V_\gamma,
%\end{cases}
%\end{align}
%and
%\begin{align}
%\begin{cases}
% \bC[x_2]/(x_2^2) & \bC x_2 \subset V_\gamma, \\
% \bC & \bC x_2 \not \subset V_\gamma.
%\end{cases}
%\end{align}
If we write an element of
$
 \Jac_{\w_\gamma} \otimes \Lambda^{\dim N_\gamma} N_\gamma^\dual
$
as
\begin{align}
 x_0^{k_0} x_1^{k_1} x_2^{k_2}
  \otimes x_{j_1}^\dual \wedge x_{j_2}^\dual \wedge \ldots \wedge x_{j_s}^\dual,
\end{align}
where $k_i = 0$ if $\bC x_i \not \subset V_\gamma$
for $i = 0, 1, 2$,
then its degree is given by
\begin{align}
 k_0 \chi_0 + k_1 \chi_1 + k_2 \chi_2 - \chi_{j_1} - \cdots - \chi_{j_s},
\end{align}
which can be proportional to $\chi$
only if
$
 V_\gamma \cap \lb \bC x_3 \oplus \cdots \bC x_{n+1} \rb
$
is either
$\bC x_3 \oplus \cdots \oplus \bC x_{n+1}$
or $0$.
We will assume this condition
for the rest of \pref{sc:E6},
and divide the analysis into the following three cases:
\begin{itemize}
  \item $\bC x_0 \subset V_\gamma$.
  \item $\bC x_0 \not \subset V_\gamma$ and $V_\gamma \neq 0$.
  \item $V_\gamma = 0$.
\end{itemize}

\subsubsection{}

Since $t_0=1$ implies $t_2 = 1$ and $t_1 = \pm 1$,
one has the following:
\begin{itemize}
 \item
$V_\gamma$ contains $\bC x_0$ if and only if either
\begin{itemize}
 \item
$\gamma = (1, \ldots, 1)$, where $V_\gamma = V$,
 \item
$\gamma = (1,1,1,-1,\ldots,-1)$ with odd $n \ge 3$,
where $V_\gamma = \bC x_0 \oplus \bC x_1 \oplus \bC x_2$,
 \item
$\gamma = (1,-1,1,-1,\ldots,-1)$ with even $n$,
where $V_\gamma = \bC x_0 \oplus \bC x_2$.
\end{itemize}
\end{itemize}

\subsubsection{}

One has
$
 V_\gamma = V
$
if and only if $\gamma$ is the identity element.
The degree of
$
 x_0^{k_0} x_1^{k_1} x_2^{k_2} \in \Jac_\w
$
is
\begin{align}
% k_0 (\chi - \chi_1 - \cdots - \chi_n) + k_1 \chi_1 =
  k_0 \chi - (k_0 - k_1) \chi_1 - (k_0 - k_2) \chi_2 - k_0 \chi_3 - \cdots - k_0 \chi_{n+1},
\end{align}
which is proportional to $\chi$
if and only if
\begin{itemize}
  \item 
  4 divides $k_0-k_1$ and
  3 divides $k_0 - k_2$
  if $n=1$, and
  \item
  4 divides $k_0-k_1$,
  3 divides $k_0 - k_2$, and
  $k_0$ is even if $n>1$. 
\end{itemize}

Thus, for $n=1$, we must have
\begin{align}
\label{eq:e6a}
5k_0 + 3k_1 + 4k_2 = 12m 
\end{align} 
for $m \in \bN$, in which case
one has
\begin{align}
 \deg x_0^{k_0} x_1^{k_1} x_2^{k_2} = m \chi.
\end{align}
For each $m \in \bN$ such that $5\nmid m$, the equation (\ref{eq:e6a}) has a unique solution with $(k_1,k_2) \in \{0,1,2\} \times \{0,1\}$ and if $5 \mid m$, then there are precisely two contributions with $(k_1,k_2)= (0,0)$ and $(k_1,k_2)=(2,1)$ such that $(k_1,k_2,m) \in \{0,1,2\} \times  \{0,1\} \times \mathbb{N}$ except if $m=0$, then only $(k_1,k_2) = (0,0)$ contributes. Each such $(k_1,k_2,m)$ contributes $
 \bC(k_0n)
$
to $\HH^{t}$ and $\HH^{t+1}$
for $t = 2m$.

For $n>1$,
the condition that $k_0$ is even
forces $k_1 \ne 1$,
and the possible $(k_0, k_1, k_2)$ and
$
 t = 2 \deg(x_0^{k_0} x_1^{k_1} x_2^{k_2})/\chi
$
are given by
\begin{align} \label{eq:kE6}
\begin{array}{ccc}
  \toprule
 (k_1, k_2) & k_0 & t \\
  \midrule
 (0,0) & 12m & 22m-12mn \\
 (0,1) & 4+12m & 8+22m-(4+12m)n \\
 (2,0) & 6+12m & 12+22m-(6+12m)n \\
 (2,1) & 10+12m & 20+22m-(10+12m)n \\
  \bottomrule
\end{array}
\end{align}
for $m \in \bN$.
Each $(k_0, k_1, k_2)$
%$x_0^{k_0} x_1^{k_1} x_2^{k_2}$
from \pref{eq:kE6} contributes
$
 \bC(k_0n)
$
to $\HH^t$ and $\HH^{t+1}$.

\subsubsection{}

One has
$
 V_\gamma = \bC x_0 \oplus \bC x_1 \oplus \bC x_2 \subsetneq V
$
for
$
\gamma \in \ker \chi
$
if and only if $n$ is an odd integer greater than or equal to $3$ and
$
 \gamma = (1, 1, 1, -1, \ldots, -1).
$
The degree of
\begin{align}
 x_0^{k_0} x_1^{k_1} x_2^{k_2} \otimes x_3^\dual \wedge \cdots \wedge x_{n+1}^\dual
  \in  \Jac_{\w_\gamma} \otimes \Lambda^{\dim N_\gamma} N_\gamma^\dual
\end{align}
is given by
\begin{align}
 k_0 \chi - (k_0 - k_1) \chi_1 - (k_0 - k_2 ) \chi_2 - (k_0+1) \chi_3 - \cdots - (k_0+1) \chi_{n+1},
\end{align}
which is proportional to $\chi$
if and only if $k_0$ is odd,
$4$ divides $k_0 - k_1$, and
$3$ divides $k_0 - k_2$.
This forces $k_1 = 1$ and the possible $(k_0, k_1, k_2)$ and
\begin{align}
 t = 2 \deg(x_0^{k_0} x_1^{k_1} x_2^{k_2} \otimes x_3^\dual \wedge \cdots \wedge x_{n+1}^\dual)/\chi + \dim N_\gamma
\end{align}
are given by
\begin{align} \label{eq:kE6_2}
\begin{array}{ccc}
  \toprule
 (k_1, k_2) & k_0 &  t \\
  \midrule
 (1,0) & 9+12m & 17+22m-(9+12m)n \\
 (1,1) & 1+12m & 3+22m-(1+12m)n \\
  \bottomrule
\end{array}
\end{align}
for $m \in \bN$.
Each $(k_0, k_1, k_2)$
%$x_0^{k_0} x_1^{k_1} x_2^{k_2}$
from \pref{eq:kE6_2} contributes
$
 \bC(k_0n)
$
to $\HH^t$ and $\HH^{t+1}$.

\subsubsection{}

One has
$
 V_\gamma = \bC x_0 \oplus \bC x_2
$
if and only if
$n$ is even and
$
 \gamma = (1, -1, 1, -1, \ldots, -1) \in \ker \chi.
$
The degree of
\begin{align}
 x_0^{k_0} x_2^{k_2} \otimes x_1^\dual \wedge x_3^\dual \wedge \cdots \wedge x_{n+1}^\dual
  \in  \Jac_{\w_\gamma} \otimes \Lambda^{\dim N_\gamma} N_\gamma^\dual
\end{align}
is given by
\begin{align}
 k_0 \chi - (k_0+1) \chi_1 - (k_0 - k_2 ) \chi_2 - (k_0+1) \chi_3 - \cdots - (k_0+1) \chi_{n+1},
\end{align}
which is proportional to $\chi$
if and only if
%$k_0$ is odd,
$4$ divides $k_0+1$ and
$3$ divides $k_0 - k_2$.
The possible $(k_0, k_2)$ and
\begin{align}
 t = 2 \deg(x_0^{k_0} x_2^{k_2} \otimes x_1^\dual \wedge x_3^\dual \wedge \cdots \wedge x_n^\dual)/\chi + \dim N_\gamma
\end{align}
are given by
\begin{align} \label{eq:kE6_3}
\begin{array}{ccc}
  \toprule
 k_2 & k_0 & t \\
  \midrule
 0 & 3+12m & 6+22m-(3+12m)n \\
 1 & 7+12m & 14+22m-(7+12m)n \\
  \bottomrule
\end{array}
\end{align}
for $m \in \bN$.
Each $(k_0, k_2)$
%$x_0^{k_0} x_1^{k_1} x_2^{k_2}$
from \pref{eq:kE6_3} contributes
$
 \bC(k_0n)
$
to $\HH^t$ and $\HH^{t+1}$.

\subsubsection{}

If $V_\gamma = \bC x_1$,
then one has
\begin{align}
 \deg \lb x_1^{k_1} \otimes x_0^\dual \wedge x_2^\dual \wedge \cdots \wedge x_{n+1}^\dual \rb
 &= - \chi_0 + k_1 \chi_1 - \chi_2 - \cdots - \chi_{n+1} \\
 &= - \chi + (k_1+1) \chi_1,
\end{align}
which is not proportional to $\chi$
for any $k_1 \in \{ 0, 1, 2 \}$.
Similarly,
$\gamma$ with $\bC x_0 \not \subset V_\gamma$
and $V_\gamma \ne 0$
does not contribute to $\HH^*$.

\subsubsection{}

One has $V_\gamma = 0$
if and only if
$
 t_1 \in \lb \bsmu_4 \setminus \{ 1 \} \rb,
$
$
 t_2 \in \lb \bmu_3 \setminus \{ 1 \} \rb,
$
and
$
 t_3 = \cdots, t_{n+1} = -1,
$
since $t_2 \ne 1$ implies
$
 t_0 = (-1)^{n-1} t_1^{-1} t_2^{-1} \ne 1.
$
There are six such $\gamma$,
and each of them contributes $\bC(-n)$
to $\HH^{n}$.

\subsection{Type $E_7$}
 \label{sc:E7}

Consider the case
\begin{align}
 \w = x_1^3 x_2 + x_2^3 + x_3^2 + \cdots + x_{n+1}^2 \in \bC[x_0, x_1,\ldots,x_{n+1}]
\end{align}
with
\begin{align}
 \Gamma
  = \Gamma_\w
  \coloneqq \lc \gamma = (t_0, \ldots, t_{n+1}) \in (\Gm)^{n+2} \relmid
   t_1^3 t_2 = t_2^3 = t_3^2 = \cdots = t_{n+1}^2 = t_0 \cdots t_{n+1} \rc,
\end{align}
so that
$
 \ker \chi \cong \bmu_9 \times \lb \bmu_2 \rb^{n-1}
$
and $\Char(\Gamma)$ is generated by
$\chi$ and $\chi_i = \deg x_i$ for $i \in \{ 0, \ldots, n+1 \}$
with relations
\begin{align}
 \chi = 3 \chi_1 + \chi_2 = 3 \chi_2 = 2 \chi_3 = \cdots = 2 \chi_{n+1} = \chi_0 + \cdots + \chi_{n+1}.
\end{align}
These relations imply
\begin{align}
 \chi_2 &= \chi - 3 \chi_1, \\
 9 \chi_1 &= 2 \chi, \\
 \chi_0
  &= \chi - \chi_1 - \cdots - \chi_{n+1} \\
  &= 2 \chi_1 - \chi_3 - \cdots - \chi_{n+1}.
\end{align}

\subsubsection{}

For any $\gamma \in \ker \chi$,
the intersection 
$
 V_\gamma \cap \lb \bC x_1 \oplus \bC x_2 \rb
$
can be either 
$\bC x_1 \oplus \bC x_2$,
$\bC x_2$, or
$0$,
where $\Jac_{\w_\gamma'}$ is isomorphic to
$
 \bC[x_1,x_2]/(3 x_1^2 x_2, x_1^3 + 3 x_2^2),
$
$
 \bC[x_2]/(3 x_2^2),
$
or
$
 \bC
$
respectively.
A basis of
$
 \bC[x_1,x_2]/(3 x_1^2 x_2, x_1^3 + 3 x_2^2)
$
is given by $\{ 1, x_1, x_1^2, x_1^3, x_1^4, x_2, x_1 x_2 \}$.

If we write an element of
$
 \Jac_{\w_\gamma} \otimes \Lambda^{\dim N_\gamma} N_\gamma^\dual
$
as
\begin{align}
 x_0^{k_0} x_1^{k_1} x_2^{k_2}
  \otimes x_{j_1}^\dual \wedge x_{j_2}^\dual \wedge \ldots \wedge x_{j_s}^\dual,
\end{align}
then its degree is given by
\begin{align}
 k_0 \chi_0 + k_1 \chi_1 + k_2 \chi_2 - \chi_{j_1} - \cdots - \chi_{j_s},
\end{align}
which can be proportional to $\chi$
only if
$
 V \cap \lb \bC x_3 \oplus \cdots \oplus \bC x_{n+1} \rb
$
is either
$
 \bC x_3 \oplus \cdots \oplus \bC x_{n+1}
$
or $0$.
We assume this condition
for the rest of \pref{sc:E7}.

\subsubsection{}

For $\gamma = (t_0,\ldots, t_{n+1}) \in \ker \chi$,
one has $t_1^2 = t_0 t_3 \cdots t_{n+1} = \pm t_0$ and
$t_2^2 = t_0 t_1 t_3 \cdots t_{n+1} = \pm t_0 t_1$,
so that the condition $t_0 = 1$ implies $t_1^2 = \pm 1$ and $t_2^2 = \pm t_1$,
which together with $t_2^3=1$ imply $t_1=t_2=1$.
Hence one has $\bC x_0 \subset V_\gamma$
if and only if either $V_\gamma = V$
or $V_\gamma = \bC x_0 \oplus \bC x_1 \oplus \bC x_2$.

\subsubsection{}

One has
$
 V_\gamma = V
$
if and only if $\gamma$ is the identity element.
The degree of
$
 x_0^{k_0} x_1^{k_1} x_2^{k_2} \in \Jac_\w
$
is
\begin{align}
 & k_0 (2 \chi_1 - \chi_3 - \cdots - \chi_{n+1}) + k_1 \chi_1 + k_2 (\chi - 3 \chi_1) \\
  &\qquad =
 k_2 \chi + (2 k_0 + k_1 - 3 k_2) \chi_1 - k_0 \chi_3 - \cdots - k_0 \chi_{n+1},
\end{align}
which is proportional to $\chi$
if and only if
\begin{itemize}
  \item 9 divides $2k_0 + k_1- 3k_2$ if $n=1$, and
  \item 9 divides $2k_0 + k_1- 3k_2$ and $k_0$ is even if $n > 1$.
\end{itemize}
For $n=1$,
one has
\begin{align}
  t
  &\coloneqq 2 \deg (x_0^{k_0} x_1^{k_1} x_2^{k_2})/\chi \\
  &= 2 k_2 + \frac{4}{9}(2k_0+k_1-3k_2) .
\end{align}
The possible $(k_0, k_1, k_2)$ and $t$
are given by
\begin{align} \label{eq:kE7_1}
\begin{array}{ccc}
  \toprule
 (k_1, k_2) & k_0 & t \\
  \midrule
 (0,0) & 9m & 8m \\
 (1,0) & 4+9m & 4+8m \\
 (2,0) & 8+9m & 8+8m \\
 (3,0) & 3+9m & 4+8m \\
 (4,0) & 7+9m & 8+8m \\
 (0,1) & 6+9m & 6+8m \\
 (1,1) & 1+9m & 2+8m \\
  \bottomrule
\end{array}
\end{align}
for $m \in \bN$.
Each $(k_0, k_1, k_2)$
from \pref{eq:kE7_1} contributes
$
 \bC(k_0n)
$
to $\HH^t$ and $\HH^{t+1}$.

In addition,
for the case $(k_1,k_2)= (2,0)$,
the element
$
 x_0^\dual \otimes x_1^2
$
corresponding to $m=-1$ in \pref{eq:kE7_1}
has degree 0,
and contributes $\bC(-1)$ to $\HH^1$.

For $n>1$,
one has
\begin{align}
 t
 &\coloneqq 2 \deg (x_0^{k_0} x_1^{k_1} x_2^{k_2})/\chi \\
 &= 2 k_2 + \frac{4}{9}(2 k_0 + k_1 - 3 k_2) - k_0 (n-1).
\end{align}
The possible $(k_0, k_1, k_2)$ and $t$ are given by
\begin{align} \label{eq:kE7}
\begin{array}{ccc}
  \toprule
 (k_1, k_2) & k_0 & t \\
  \midrule
 (0,0) & 18m & 34m-18mn \\
 (1,0) & 4+18m & 8+34m-(4+18m)n \\
 (2,0) & 8+18m & 16+34m-(8+18m)n \\
 (3,0) & 12+18m & 24+34m-(12+18m)n \\
 (4,0) & 16+18m & 32+34m-(16+18m)n \\
 (0,1) & 6+18m & 12+34m-(6+18m)n \\
 (1,1) & 10+18m & 20+34m-(10+18m)n \\
  \bottomrule
\end{array}
\end{align}
for $m \in \bN$.
Each $(k_0, k_1, k_2)$
%$x_0^{k_0} x_1^{k_1} x_2^{k_2}$
from \pref{eq:kE7} contributes
$
 \bC(k_0n)
$
to $\HH^t$ and $\HH^{t+1}$.

\subsubsection{}

For $n>1$, in addition, one has
$
 V_\gamma = \bC x_0 \oplus \bC x_1 \oplus \bC x_2
$
if and only if $n$ is odd and
$\gamma = (1,1,1,-1,\ldots,-1)$.
The degree of
\begin{align}
 x_0^{k_0} x_1^{k_1} x_2^{k_2} \otimes x_3^\dual \wedge \cdots \wedge x_{n+1}^\dual
 \in \Jac_{\w_\gamma} \otimes \Lambda^{\dim N_\gamma} N_\gamma^\dual
\end{align}
is given by
\begin{align}
 & k_0 (2 \chi_1 - \chi_3 - \cdots - \chi_{n+1}) + k_1 \chi_1 + k_2 (\chi - 3 \chi_1)
  - \chi_3 - \cdots - \chi_{n+1} \\
  &\qquad =
 k_2 \chi + (2 k_0 + k_1 - 3 k_2) \chi_1 - (k_0+1) \chi_3 - \cdots - (k_0+1) \chi_{n+1},
\end{align}
which is proportional to $\chi$
if and only if
9 divides $2 k_0 + k_1 - 3 k_2$ and
$k_0$ is odd.
% in which case the constant of proportionality is
% \begin{align}
%  k_2 + \frac{2}{9}(2 k_0 + k_1 - 3 k_2) - \frac{1}{2} (k_0+1) (n-1).
% \end{align}
The possible $(k_0, k_1,k_2)$ and
\begin{align}
 t \coloneqq 2 \deg \lb x_0^{k_0} x_1^{k_1} x_2^{k_2} \otimes x_3^\dual \wedge \cdots \wedge x_{n+1}^\dual \rb / \chi + \dim N_\gamma
\end{align}
are given by
\begin{align} \label{eq:kE7_2}
\begin{array}{ccc}
  \toprule
 (k_1, k_2) & k_0 & t \\
  \midrule
 (0,0) & 9+18m & 17+34m-(9+18m)n\\
 (1,0) & 13+18m & 25+34m-(13+18m)n\\
 (2,0) & 17+18m & 33+34m-(17+18m)n\\
 (3,0) & 3+18m & 7+34m-(3+18m)n\\
 (4,0) & 7+18m & 15+34m-(7+18m)n\\
 (0,1) & 15+18m & 29+34m-(15+18m)n\\
 (1,1) & 1+18m & 3+34m-(1+18m)n\\
  \bottomrule
\end{array}
\end{align}
for $m \in \bN$.
Each $(k_0, k_1, k_2)$
%$x_0^{k_0} x_1^{k_1} x_2^{k_2}$
from \pref{eq:kE7_2} contributes
$
 \bC(k_0n)
$
to $\HH^t$ and $\HH^{t+1}$.

In addition,
for the case $(k_1,k_2)= (2,0)$,
the element
$
 x_0^\dual \otimes x_1^2 \otimes
  x_3^\dual \wedge \cdots \wedge x_{n+1}^\dual
$
corresponding to $m=-1$ in \pref{eq:kE7_2}
has degree 0,
and contributes $\bC(-n)$ to $\HH^n$.

\subsubsection{}

One has $V_\gamma=0$
for $\gamma = (t_0, \ldots, t_{n+1}) \in \ker \chi$
if and only if
$
 t_1 \in \bmu_9 \setminus \{ 1 \},
$
$
 t_2 \coloneqq t_1^{-3} \ne 1,
$
and
$
 t_3 = \cdots = t_{n+1} = -1,
$
in which case one has
$
 t_0 = (-1)^{n-1} t_1^2 \ne 1.
$
The set $\lc t_1 \in \bmu_9 \relmid t_1^3 \ne 1 \rc$
consists of six elements,
each of which contributes
$\bC(-n)$
to $\HH^{n}$.

\subsection{Type $E_8$}
 \label{sc:E8}

Consider the case
\begin{align}
 \w = x_1^5 + x_2^3 + x_3^2 + \cdots + x_{n+1}^2 \in \bC[x_0, x_1,\ldots,x_{n+1}]
\end{align}
with
\begin{align}
 \Gamma
  = \Gamma_\w
  \coloneqq \lc \gamma = (t_0, \ldots, t_{n+1}) \in (\Gm)^{n+2} \relmid
   t_1^5 = t_2^3 = t_3^2 = \cdots = t_{n+1}^2 = t_0 \cdots t_{n+1} \rc,
\end{align}
so that
$
 \ker \chi \cong \bmu_5 \times \bmu_3 \times \lb \bmu_2 \rb^{n-1}
$
and $\Char(\Gamma)$ is generated by
$\chi$ and $\chi_i = \deg x_i$ for $i \in \{ 0, \ldots, n+1 \}$
with relations
\begin{align}
 \chi = 5 \chi_1 = 3 \chi_2 = 2 \chi_3 = \cdots = 2 \chi_{n+1} = \chi_0 + \cdots + \chi_{n+1}.
\end{align}

\subsubsection{}

If we write an element of
$
 \Jac_{\w_\gamma} \otimes \Lambda^{\dim N_\gamma} N_\gamma^\dual
$
as
\begin{align}
 x_0^{k_0} x_1^{k_1} x_2^{k_2}
  \otimes x_{j_1}^\dual \wedge x_{j_2}^\dual \wedge \ldots \wedge x_{j_s}^\dual,
\end{align}
then its degree is given by
\begin{align}
 k_0 \chi_0 + k_1 \chi_1 + k_2 \chi_2 - \chi_{j_1} - \cdots - \chi_{j_s},
\end{align}
which can be proportional to $\chi$
only if
$
 V \cap \lb \bC x_3 \oplus \cdots \oplus \bC x_{n+1} \rb
$
is either
$
 \bC x_3 \oplus \cdots \oplus \bC x_{n+1}
$
or $0$.
We assume this condition
for the rest of \pref{sc:E8}.

\subsubsection{}

Since $t_0 = 1$ implies $t_1=t_2=1$,
one has $\bC x_0 \subset V_\gamma$
if and only if either $V_\gamma = V$
or $V_\gamma = \bC x_0 \oplus \bC x_1 \oplus \bC x_2$.

\subsubsection{}

One has
$
 V_\gamma = V
$
if and only if $\gamma$ is the identity element.
The degree of
$
 x_0^{k_0} x_1^{k_1} x_2^{k_2} \in \Jac_\w
$
is
\begin{align}
 k_0 \chi - (k_0-k_1) \chi_1 - (k_0-k_2) \chi_2 - k_0 \chi_3 - \cdots - k_0 \chi_{n+1},
\end{align}
which is proportional to $\chi$
if and only if
\begin{itemize}
  \item $5$ divides $k_0-k_1$ and $3$ divides $k_0-k_2$ if $n=1$, and
  \item $5$ divides $k_0-k_1$, $3$ divides $k_0-k_2$, and $k_0$ is even if $n>1$.
\end{itemize}
% \begin{align}
% k_0 &\equiv k_1 \mod 5, \\
% k_0 &\equiv k_2 \mod 3, 
% \end{align}
% if $n=1$, and
% \begin{align}
% k_0 &\equiv k_1 \mod 5, \\
% k_0 &\equiv k_2 \mod 3, \\
% k_0 &\equiv 0 \mod 2,
% \end{align}
% if $n>1$.

For $n=1$, we must have
\begin{align}\label{eq:e8a} 
% k_0 = 6 k_1 + 10 k_2 + 30m
7k_0 + 3 k_1 + 5k_2 = 15 m \end{align}
for $m \in \bN$, 
in which case one has
\begin{align}
%t/2 = 6k_1+10k_2+29m-(3k_1+5k_2+15m)n.
 t \coloneqq 2 \deg \lb x_0^{k_0} x_1^{k_1} x_2^{k_2} \rb/\chi = 2 m .
\end{align}
For each $m \in \mathbb{N}$ such that $7 \nmid m$, the equation (\ref{eq:e8a}) has a unique solution with $(k_1,k_2) \in \{0,1,2,3\} \times \{0,1\}$ and 
if $7 \mid m$, then there are precisely two contributions with $(k_1,k_2)=(0,0)$ and $(k_1,k_2)= (3,1)$ such that $(k_1, k_2, m) \in \{ 0, 1, 2, 3 \} \times \{ 0, 1 \} \times \bN$ except if $m=0$, then only $(k_1,k_2)=(0,0)$ contributes. 

For $n>1$, we must have
\begin{align}\label{eq:e8b} 7k_0 + 3 k_1 + 5 k_2 = 15m \end{align}
for $m \in \bN$,
and in addition $k_0$ must be in $2\bN$. Thus, we can re-write (\ref{eq:e8b}) as
\begin{align} \label{eq:e8b2}
  k_0 = 6k_1 + 10k_2+ 30m'
\end{align} 
with $m' = k_0/2 - m$.
One has
\begin{align}
  t &\coloneqq 2 \deg \lb x_0^{k_0} x_1^{k_1} x_2^{k_2} \rb/\chi \\
  % &= 2 \lb
  %  (6k_1 + 10k_2+ 30m')
  %   - \frac{1}{5}(5k_1+10k_2+30m')
  %   - \frac{1}{3}(6k_1 + 9k_2+ 30m')
  %   - \frac{1}{2}(n-1)(6k_1 + 10k_2+ 30m')
  %   \rb \\
  % &= \lb
  %  (12 k_1 + 20 k_2+ 60 m')
  %   - (2k_1+4k_2+12m')
  %   - (4k_1 + 6k_2+ 20m')
  %   - (n-1)(6k_1 + 10k_2+ 30m')
  %   \rb \\
  &= 12k_1+20k_2 +58m' - (6k_1+10k_2 +30m')n.
\end{align}
Each $(k_1,k_2,m') \in \{0,1,2,3\} \times \{0,1\} \times \bN$ contributes
$
 \bC(k_0n)
$
to $\HH^t$ and $\HH^{t+1}$.

\subsubsection{}

If $n>1$, in addition, one has
$
 V_\gamma = \bC x_0 \oplus \bC x_1 \oplus \bC x_2
$
if and only if $n$ is odd and
$\gamma = (1,1,1,-1,\ldots,-1)$.
The degree of
\begin{align}
 x_0^{k_0} x_1^{k_1} x_2^{k_2} \otimes x_3^\dual \wedge \cdots \wedge x_{n+1}^\dual
 \in \Jac_{\w_\gamma} \otimes \Lambda^{\dim N_\gamma} N_\gamma^\dual
\end{align}
is
\begin{align}
 k_0 \chi - (k_0-k_1) \chi_1 - (k_0-k_2) \chi_2 - (k_0+1) \chi_3 - \cdots - (k_0+1) \chi_{n+1},
\end{align}
which is proportional to $\chi$
if and only if
%\begin{align}
% k_0 &\equiv k_1 \mod 5, &
% k_0 &\equiv k_2 \mod 3, &
% k_0 &\equiv 1 \mod 2.
%\end{align}
%which implies
\begin{align}
\label{eq:e8c}
14k_0 + 6k_1 + 10k_2 = 30 m 
%k_0 = 15 + 6 k_1 + 10 k_2 + 30m
\end{align}
for $m \in \bZ$ and in addition we must have $k_0$ odd.
Thus, again we can rewrite \pref{eq:e8c} as 
\begin{align}
\label{eq:e8c2}
k_0 = 15 + 6 k_1 + 10 k_2 + 30m' 
\end{align}
where $m' = (k_0 -1)/2 -m$.
One has
\begin{align}
  t
  &\coloneqq
  2 \deg \lb x_0^{k_0} x_1^{k_1} x_2^{k_2}
   \otimes x_3^\dual \wedge \cdots \wedge x_{n+1}^\dual \rb / \chi
   + \dim N_\gamma \\
  &= 2 \lb
  k_0 - \frac{1}{5}(k_0-k_1) - \frac{1}{3}(k_0-k_2) - \frac{1}{2}(k_0+1)(n-1) \rb + (n-1) \\
  &= 29 + 12k_1 + 20k_2 + 58m' - (15+6k_1+10k_2+30m')n
\end{align}
Each $(k_1,k_2,m') \in \{0,1,2,3\} \times \{0,1\} \times \bZ$ such that
\begin{align}
  15 + 6k_1+ 10 k_2 + 30m ' \geq 0
\end{align}
contributes
$
 \bC(k_0n)
$
to $\HH^t$ and $\HH^{t+1}$.

%Each $(k_1, k_2, m) \in \{ 0, 1, 2, 3 \} \times \{ 0, 1 \} \times \bN$
%contributes
%$
% \bC(k_0n)
%$
%to $\HH^t$ and $\HH^{t+1}$.

\subsubsection{}

An element $\gamma = (t_0, \ldots, t_{n+1}) \in \ker \chi$
satisfies $V_\gamma = 0$
if and only if
$
 t_1 \in \bmu_5 \setminus \{ 1 \},
$
$
 t_2 \in \bmu_3 \setminus \{ 1 \},
$
$
 t_3 = \cdots = t_{n+1} = -1,
$
and
$
 t_0 = (-1)^{n-1} (t_1 t_2)^{-1}.
$
There are eight such elements,
each of which contributes
$\bC(-n)$
to $\HH^{n}$.

\bibliographystyle{amsalpha}
\bibliography{bibs}

\end{document}